\documentclass{article}
\usepackage{fullpage}
\usepackage{fancyhdr}
\usepackage{amsmath, amssymb, latexsym, amscd, amsthm,amsfonts,amstext}
\usepackage{subfigure}
\usepackage{enumerate}
\usepackage{xcolor}
\usepackage{textcomp}
\usepackage{mathtools}

\usepackage{graphicx}
\newtheorem{theorem}{Theorem}[section]
\newtheorem{lemma}[theorem]{Lemma}
\newtheorem{proposition}[theorem]{Proposition}
\newtheorem{corollary}[theorem]{Corollary}
\numberwithin{equation}{section}
\newtheorem{definition}[theorem]{Definition}
\newtheorem{remark}[theorem]{Remark}

\numberwithin{equation}{section}

\usepackage[mathscr]{eucal}
\usepackage{graphicx}
\usepackage{subfigure}
\usepackage{cite}

\usepackage{xcolor}

\usepackage{url}

\allowdisplaybreaks
\title{Multiscale analysis of the acoustic scattering by many scatterers of impedance type} 

\author{
 Durga Prasad Challa%\footnote{C\lowercase{orresponding author}: D\lowercase{urga} P\lowercase{rasad} C\lowercase{halla}.} 
\thanks{Department of Mathematics, 
Tallinn University of Technology,  Tallinn, Estonia.
(Email: durga.challa@ttu.ee).
}%\footnote{Corresponding author: Durga Prasad Challa} 
\and  Mourad Sini
\thanks{Radon institute (RICAM), Austrian Academy of Sciences, 69 Altenbergerstrasse, A4040, Linz, Austria
(Email: mourad.sini@oeaw.ac.at).}
}

\begin{document}
\graphicspath{{Figures-eps/}}
\maketitle

\begin{abstract}
 We are concerned with the acoustic scattering problem, at a frequency $\kappa$, by many small obstacles of arbitrary shapes with impedance boundary condition. 
 These scatterers are assumed to be included in a bounded domain $\Omega$
 in $\mathbb{R}^3$ which is embedded in an acoustic background characterized by an eventually locally varying index of refraction. 
 The collection of the scatterers $D_m, \; m=1,...,M$ is modeled by four parameters: their number $M$, their maximum radius $a$, their minimum distance $d$ and the surface impedances 
 $\lambda_m, \; m=1,...,M$. We consider the parameters $M, d$ and $\lambda_m$'s having the following scaling properties: $M:=M(a)=O(a^{-s})$, $d:=d(a)\approx a^t$ and 
 $\lambda_m:=\lambda_m(a)=\lambda_{m,0}a^{-\beta}$, as $a \rightarrow 0$, with non negative constants $s, t$ and $\beta$ and complex numbers $\lambda_{m, 0}$'s 
 with eventually negative imaginary parts.
 
 We derive the asymptotic expansion of the farfields with explicit error estimate in terms of $a$, as $a\rightarrow 0$. The dominant term is the Foldy-Lax field corresponding to the scattering
 by the point-like scatterers located at the centers $z_m$'s of the scatterers $D_m$'s with $\lambda_m \vert \partial D_m\vert$ as the related scattering coefficients. 
 This asymptotic expansion is justified under the following conditions $$a \leq a_0,\; \vert \Re (\lambda_{m,0})\vert \geq \lambda_-,\; \vert \lambda_{m,0}\vert \leq \lambda_+,\;  
 ~~ \beta<1, \;~~ 0\leq s \leq  2-\beta,\;~~\frac{s}{3}\leq t$$ and the error of the approximation is $C\; a^{3-2\beta-s}$, as $a \rightarrow 0$, 
 where the positive constants $a_0,\; \lambda_-,\; \lambda_+$ and $C$ depend only on the a priori uniform bounds of the Lipschitz characters of the obstacles $D_m$'s and the ones of
  $M(a)a^s$ and $\frac{d(a)}{a^t}$. We do not assume the periodicity in distributing the small scatterers. In addition,  the scatterers can be arbitrary close since $t$ can be arbitrary large, i.e. we can handle the mesoscale regime. Finally, for spherical scatterers, we can also allow the limit case $\beta=1$ with a slightly better error of the approximation. 
\end{abstract}

%\begin{keywords}

%\begin{AMS}
 %35J08, 35Q61, 45Q05 
%\end{AMS}

\pagestyle{myheadings}
 \thispagestyle{plain}
%  \markboth{D. P. Challa and M. Sini }{Acoustic scattering by scatterers of impedance type}
%\markboth{D. P. Challa and M. Sini }{Acoustic scattering by impedance type scatterers}
\section{Introduction and statement of the results}\label{Introduction-smallac-sdlp}

 Let $B_1, B_2,\dots, B_M$ be $M$ open, bounded and simply connected sets in $\mathbb{R}^3$ with Lipschitz boundaries
containing the origin. We assume that the Lipschitz constants of $B_j$, $j=1,..., M$ are uniformly bounded.  
We set $D_m:=\epsilon B_m+z_m$ to be the small bodies characterized by the parameter 
$\epsilon>0$ and the locations $z_m\in \mathbb{R}^3$, $m=1,\dots,M$. Let $U^{i}$ be a solution of the Helmholtz equation $(\Delta + \kappa^{2})U^{i}=0 \mbox{ in } \mathbb{R}^{3}$.  
We denote by  $U^{s}$ the acoustic field scattered by the $M$ small bodies $D_m\subset \mathbb{R}^{3}$, due to 
the incident field $U^{i}$ (mainly the plane incident waves $U^{i}(x,\theta):=e^{ikx\cdot\theta}$
with the incident direction $\theta \in \mathbb{S}^2$, where $\mathbb{S}^2$ being the unit sphere), with impedance boundary conditions. Hence the total field $U^{t}:=U^{i}+U^{s}$ satisfies the 
following exterior impedance problem of the acoustic waves
\begin{equation}
(\Delta + \kappa^{2})U^{t}=0 \mbox{ in }\mathbb{R}^{3}\backslash \left(\mathop{\cup}_{m=1}^M \bar{D}_m\right),\label{acimpoenetrable}
\end{equation}
\begin{equation}
\left.\frac{\partial U^{t}}{\partial \nu_m}+\lambda_m U^{t}\right|_{\partial D_m}=0,\, 1\leq m \leq M, \label{acgoverningsupport}  
\end{equation}
\begin{equation}
\frac{\partial U^{s}}{\partial |x|}-i\kappa U^{s}=o\left(\frac{1}{|x|}\right), |x|\rightarrow\infty, ~(\text{S.R.C}) \label{radiationc}
\end{equation}
where $\nu_m$ is the outward unit normal vector of $\partial D_m$, $\kappa>0$ is the wave number and S.R.C stands for the Sommerfeld radiation condition.  The scattering problem 
(\ref{acimpoenetrable}-\ref{radiationc}) is well posed in the H\"{o}lder or Sobolev spaces, see \cite{C-K:1983, C-K:1998, Mclean:2000, Ned:Spring2001} for instance,
 in the case $\Im \lambda_m \geq 0$. We will see that this condition can be relaxed to 
 allow $\Im \lambda_m$ to be negative under some conditions. Applying Green's formula to $U^s$, we can show that
 the scattered field $U^s(x, \theta)$ has the following asymptotic expansion:
\begin{equation}\label{far-field}
 U^s(x, \theta)=\frac{e^{i \kappa |x|}}{4\pi|x|}U^{\infty}(\hat{x}, \theta) + O(|x|^{-2}), \quad |x|
\rightarrow \infty,
\end{equation}
with $\hat{x}:=\frac{x}{\vert x\vert}$, where the function
$U^{\infty}(\hat{x}, \theta)$ for $(\hat{x}, \theta)\in \mathbb{S}^{2} \times\mathbb{S}^{2}$  is called the far-field pattern.
We recall that the fundamental solution, $\Phi_\kappa(x,y)$, of the Helmholtz equation in $\mathbb{R}^3$
with the fixed wave number $\kappa$ is given by $
 \Phi_\kappa(x,y):=\frac{e^{i\kappa|x-y|}}{4\pi|x-y|},\quad \text{for all } x,y\in\mathbb{R}^3.$
 
\begin{definition} %[Definition \ref{Introduction-smallac-sdlp}.\ref{Def1}.] 
\label{Def1}
We define 
\begin{enumerate}%[\ref{Def1}.1.]
 \item $
a:=\max\limits_{1\leq m\leq M } diam (D_m) ~~\big[=\epsilon \max\limits_{1\leq m\leq M } diam (B_m)\big],$

 \item  $
d:=\min\limits_{\substack{m\neq j\\1\leq m,j\leq M }} d_{mj},$
where $\,d_{mj}:=dist(D_m, D_j)$. 

\item $\kappa_{\max}$ as the upper bound of the used wave numbers, i.e. $\kappa\in[0,\,\kappa_{\max}]$.
\bigskip

The distribution of the scatterers is modeled as follows:
\item the number $M~:=~M(a)~:=~O(a^{-s})\leq M_{max} a^{-s}$ with a given positive constant $M_{max}$\\
\item the minimum distance $d~:=~d(a)~\approx ~a^t$, i.e. $d_{min} a^t \leq d(a) \leq d_{max}a^t $, with given positive constants $d_{min}$ and $d_{max}$. \\
 \item the surface impedance $\lambda_m~:=~\lambda_{m,0}a^{-\beta}$, where $\lambda_{m,0} \neq 0$ and might be a complex number.
%  such $\vert \lambda_{m, 0}\vert \leq \lambda_0$
%  with a given constant $\lambda_0$.\\
\end{enumerate}

\end{definition}
Here the real numbers $s$, $t$ and $\beta$ are assumed to be non negative.
\bigskip

We call the upper bounds of the Lipschitz character of $B_m$'s, $M_{max}, d_{min}, d_{max}$ 
and $\kappa_{max}$ the set of the a priori bounds. 
\par 
The goal of our work is to derive an asymptotic expansion of the scattered field by the collection of the small scatterers $D_m,m=1,\dots,M$, taking into account these
parameters. This is the object of the following theorem.
\bigskip

\begin{theorem}\label{Maintheorem-ac-small-sing}
 There exist positive constants $a_0$, $\lambda_-$ and $\lambda_+$ depending only on the set
 of the a priori bounds such that if 
\begin{equation} \label{conditions}
a \leq a_0,\; \vert \lambda_{m, 0}\vert \leq \lambda_+, \; \vert \Re(\lambda_{m, 0}) \vert \geq \lambda_-,\; ~~ \beta<1, \;~~ s \leq 2-\beta,\;~~\frac{s}{3}\leq t
% \,,  \;~~ \min_{j\neq m}\cos(\kappa \vert z_j- z_m\vert)\geq 0
\end{equation} 
then the far-field pattern $U^\infty(\hat{x},\theta)$ has the following asymptotic expansion
\begin{equation}\label{x oustdie1 D_m farmain}
U^\infty(\hat{x},\theta)=\sum_{m=1}^{M}e^{-i\kappa\hat{x}\cdot z_m}Q_m+
O\left(
 a^{3-s-2\beta}\right),
 \end{equation}
uniformly in $\hat{x}$ and $\theta$ in $\mathbb{S}^2$. The constant appearing in the estimate $O(.)$ depends only on 
the set of the a priori bounds, $\lambda_-$ and $\lambda_+$. The coefficients $Q_m$, $m=1,..., M,$ are the solutions of the following linear algebraic system
\begin{eqnarray}\label{fracqcfracmain}
 Q_m +\sum_{\substack{j=1 \\ j\neq m}}^{M}C_m \Phi_\kappa(z_m,z_j)Q_j&=&-C_mU^{i}(z_m, \theta),~~
\end{eqnarray}
for $ m=1,..., M,$ with 
$C_m:=-\lambda_m|\partial D_m|$.%$\frac{\lambda|\partial D_m|}{\left[-1+\lambda J_m\right]}$ and  $J_m:=\left\{\begin{array}{ccc}
%   \int_{\partial D_m}\Phi_0(s_m,t) d{s_m},\,t\in\partial D_m,&&\mbox{if $D_m$ are balls}\\
%   0, &&\mbox{else}
%   \end{array}\right.$.
  \par
%where $\bar{Q}_m:=\int_{\partial D_m}\sigma_m(s)ds$.
The algebraic system \eqref{fracqcfracmain} is invertible under 
the conditins:%two conditions:
%\footnote{If the real parts of $\lambda_{m,\; 0}$ are negative then we can get rid of the condition $s \leq 2 -\beta$ and replace it by the following $\min_{j\neq m}\cos(\kappa \vert z_j- z_m\vert)\geq 0$.}

\begin{eqnarray}\label{invertibilityconditionsmainthm}
 %t\leq 2-\beta \text{ and } 
 \vert \Re(\lambda_{m,0})\vert \geq \lambda_-,\; s \leq  2 -\beta.
\end{eqnarray} 

\end{theorem}
 
Let us now assume that the background is not homogeneous but modeled 
by a locally variable index of refraction $n$, i.e. there exists a bounded set $\Omega$ such that $n(x)=1,\; x\in \mathbb{R}^{3}\setminus \Omega$ and bounded inside $\Omega$, i.e. $\vert n(x)\vert \leq n_{max},\; x\in \Omega$. In this case we model our scattering problem as follows. The total field $U_n^{t}:=U^{i}+U_n^{s}$ satisfies the 
following exterior impedance problem of the acoustic waves
\begin{equation}
(\Delta + \kappa^{2}n^2(x))U_n^{t}=0 \mbox{ in }\mathbb{R}^{3}\backslash \left(\mathop{\cup}_{m=1}^M \bar{D}_m\right),\label{acimpoenetrable-1}
\end{equation}
\begin{equation}
\left.\frac{\partial U_n^{t}}{\partial \nu_m}+\lambda_m U_n^{t}\right|_{\partial D_m}=0,\, 1\leq m \leq M, \label{acgoverningsupport-1}  
\end{equation}
\begin{equation}
\frac{\partial U_n^{s}}{\partial |x|}-i\kappa U_n^{s}=o\left(\frac{1}{|x|}\right), |x|\rightarrow\infty, \label{radiationc-1}
\end{equation}
Again, the scattering problem 
(\ref{acimpoenetrable-1}-\ref{radiationc-1}) is well posed in the H\"{o}lder or Sobolev spaces, see \cite{C-K:1983, C-K:1998, Mclean:2000}
 in the case $\Im \lambda_m >0$. As we said for (\ref{acimpoenetrable-1}-\ref{radiationc-1}), this last condition can be relaxed to allow $\Im \lambda_m$ to be negative. 
 Applying Green's formula to $U_n^s$, we can show that
 the scattered field $U_n^s(x, \theta)$ has the following asymptotic expansion:
\begin{equation}\label{far-field-n}
 U_n^s(x, \theta)=\frac{e^{i \kappa |x|}}{4\pi|x|}U_n^{\infty}(\hat{x}, \theta) + O(|x|^{-2}), \quad |x|
\rightarrow \infty,
\end{equation}
where the function
$U_n^{\infty}(\hat{x}, \theta)$ for $(\hat{x}, \theta)\in \mathbb{S}^{2} \times\mathbb{S}^{2}$  is the corresponding far-field pattern.
As a corollary of Theorem \ref{Maintheorem-ac-small-sing}, we derive the following result:

\begin{corollary}\label{Maintheorem-ac-small-sing-1}
 There exist positive constants $a_0$, $\lambda_-$, $\lambda_+$ depending only on the 
set of the a priori bounds and on $n_{max}$ such that
if 
\begin{equation} \label{conditions-1}
a \leq a_0,\; \vert \lambda_{m, 0}\vert \leq \lambda_+, \; \vert \Re(\lambda_{m, 0}) \vert \geq \lambda_-,\;~~ \beta<1, \;~~ s \leq  2-\beta,\;~~\frac{s}{3}\leq t
% \,,  \;~~ \min_{j\neq m}\cos(\kappa \vert z_j- z_m\vert)\geq 0
\end{equation} 
then the far-field pattern $U_n^\infty(\hat{x},\theta)$ has the following asymptotic expansion 
\footnote{ The following remarks apply for both Theorem \ref{Maintheorem-ac-small-sing} and Corollary \ref{Maintheorem-ac-small-sing-1}. 
If $s< 2-\beta$, then the bounds $\lambda_-$ and $\lambda_+$ can be arbitrary. 
If $\Re \lambda_m \leq 0$ then, for the invertibility of the systems (\ref{fracqcfracmain}) and (\ref{fracqcfracmain-1}), the condition $s\leq 2-\beta$ on the number of the small bodies can be replaced by the condition $t\leq 2-\beta$ on the minimum distance between them, see Remark \ref{Remark-on-sign-lambda}.   }
\begin{equation}\label{x oustdie1 D_m farmain-1}
U_n^\infty(\hat{x},\theta)=V_n^\infty(\hat{x},\theta)+\sum_{m=1}^{M}V_n^t(z_m, -\hat{x})\mathbb{Q}_m+
O\left(
 a^{3-s-2\beta}\right),
 \end{equation}
uniformly in $\hat{x}$ and $\theta$ in $\mathbb{S}^2$. The constant appearing in the estimate $O(.)$ depends only on 
the set of the a priori bounds, $\lambda_-$, $\lambda_+$ and $n_{max}$. The quantity\; 
\footnote{ By the mixed reciprocity relations, we have $V_n^t(z_m, -\hat{x})=G^{\infty}_{\kappa}(\hat{x}, z_m)$ where $G^{\infty}_{\kappa}(\hat{x}, z_m)$ is 
the farfield created by a point source $G_{\kappa}(x, z_m)$, located at the point $z_m$, in the scattering model (\ref{acimpoenetrable-2}-\ref{radiationc-2}).} 
$V^t(z_m, -\hat{x})$ is the total field, evaluated at the point $z_m$ in the direction $-\hat{x}$, corresponding to the scattering problem 
\begin{equation}
(\Delta + \kappa^{2}n^2(x))V_n^{t}=0 \mbox{ in }\mathbb{R}^{3},\label{acimpoenetrable-2}
\end{equation}
\begin{equation}
\frac{\partial V_n^{s}}{\partial |x|}-i\kappa V_n^{s}=o\left(\frac{1}{|x|}\right), |x|\rightarrow\infty, \label{radiationc-2}
\end{equation}
 The coefficients $\mathbb{Q}_m$, $m=1,..., M,$ are the solutions of the following linear algebraic system
\begin{eqnarray}\label{fracqcfracmain-1}
 \mathbb{Q}_m +\sum_{\substack{j=1 \\ j\neq m}}^{M}C_m G_\kappa(z_m,z_j)\mathbb{Q}_j&=&-C_mV^{t}(z_m, \theta),~~
\end{eqnarray}
for $ m=1,..., M,$ with 
\begin{equation}\label{Cm-s}
C_m:=-\lambda_m|\partial D_m|.
\end{equation} Here $G_{\kappa}(x, z)$ is the outgoing Green's function corresponding to the scattering problem (\ref{acimpoenetrable-2}-\ref{radiationc-2}).
  \par
The algebraic system \eqref{fracqcfracmain-1} is invertible under the 
conditions:
%\footnote{As in the homogeneous background case, if the real parts of $\lambda_{m,\; 0}$ are negative then we can get rid of the condition $s \leq 2 -\beta$ and replace it by the following one $\min_{j\neq m}\cos(\kappa \vert z_j- z_m\vert)\geq 0$.}

\begin{eqnarray}\label{invertibilityconditionsmainthm-1}
 %t\leq 2-\beta \text{ and } 
 \vert \Re(\lambda_{m, 0}) \vert \geq \lambda_-,\;~~ s \leq 2 -\beta.
\end{eqnarray} 

\end{corollary} 

Before discussing these results comparing them to some of  the literature, we add the following remark on the particular case when the scatterers are spherical.
\begin{remark}\label{Remark-spherical-scatterers}
If the scatterers have spherical shapes, then we can slightly improve the error estimates in (\ref{x oustdie1 D_m farmain-1}) by slightly changing the form of the algebraic system (\ref{fracqcfracmain-1}). Precisely we have $O\left(
 a^{3-s-\beta}\right)$ instead of $O\left(
 a^{3-s-2\beta}\right)$. For these spherical shapes, we can also handle the case $\beta=1$. The result reads as follows.
There exist positive constants $a_0$, $\lambda_-$, $\lambda_+$ depending only on the 
set of the a priori bounds and on $n_{max}$ such that
if 
\begin{equation} \label{conditions-1Re}
a \leq a_0,\; \vert \lambda_{m, 0}\vert \leq \lambda_+, \; \vert \Re(\lambda_{m, 0}) \vert \geq \lambda_-,\;~~ \beta \leq 1, \;~~ s \leq  2-\beta,\;~~\frac{s}{3}\leq t
% \,,  \;~~ \min_{j\neq m}\cos(\kappa \vert z_j- z_m\vert)\geq 0
\end{equation} 
then the
the far-field pattern $U_n^\infty(\hat{x},\theta)$ has the following asymptotic expansion 

\begin{equation}\label{x oustdie1 D_m farmain-2}
U_n^\infty(\hat{x},\theta)=V_n^\infty(\hat{x},\theta)+\sum_{m=1}^{M}V_n^t(z_m, -\hat{x})\mathbb{Q}_m+
O\left(
 a^{3-s-\beta}\right),
 \end{equation}
uniformly in $\hat{x}$ and $\theta$ in $\mathbb{S}^2$. The coefficients $\mathbb{Q}_m$, $m=1,..., M,$ are the solutions of the following linear algebraic system
\begin{eqnarray}\label{fracqcfracmain-2}
 \mathbb{Q}_m +\sum_{\substack{j=1 \\ j\neq m}}^{M}C_m G_\kappa(z_m,z_j)\mathbb{Q}_j&=&-C_mV^{t}(z_m, \theta),~~
\end{eqnarray}
for $ m=1,..., M,$ where now we have a slight change in the constant $C_m$ compared to (\ref{Cm-s}), i.e.
\begin{equation}\label{Cm-s-1}
C_m:=\frac{\lambda_m|\partial D_m|}{-1+\lambda_m I_m}.
\end{equation}
where $I_m:=\int_{\partial D_m}\Phi_0(s_m,t) d{s_m},\; t\in\partial D_m$, is a constant if $D_m$ is a ball.
\footnote{ The property that $I_m$ is a constant only for balls is known as Gruber's conjecture and it is solved 
for the 2D for Lipschitz regular domains and in 3D for convex domains, see \cite{M-R:2000}. 
This is the main restriction why we cannot obtain in Theorem \ref{Maintheorem-ac-small-sing} 
and Corollary \ref{Maintheorem-ac-small-sing-1} the same results 
as in Remark \ref{Remark-spherical-scatterers} for general shapes. 
Actually, if $D_m$ is a ball then the value of $I_m$ is nothing but the radius of $D_m\,(= \epsilon\;\operatorname{radius}(B_m))$.}

This algebraic system is invertible under the same conditions as 
(\ref{fracqcfracmain-1}).
\end{remark}
\bigskip

These results say that the dominant term in the approximation (\ref{x oustdie1 D_m farmain}), and similarly the one in (\ref{x oustdie1 D_m farmain-1}), 
is the Foldy-Lax field corresponding to the scattering
 by the point-like scatterers located at the 'centers' $z_m$'s of the scatterers $D_m$'s with $\lambda_m \vert \partial D_m\vert$ as the related scattering coefficients.
 This Foldy-Lax field describes the field that results in the multiple scattering between the different scatterers $D_m$'s, see \cite{Martin:2006} and the references therein 
 for more information on this issue. The accuracy of the approximation of the scattered field by the Foldy-Lax field depends of course on the error estimates. In its generality, 
 this issue is still largely open but there is an increase of interest to understand it, see for instance 
 \cite{Bendali-Coquet-Tordeux:NAA2012, Bendali-et.al:2015,DPC-SM:MMS2013, 
 DPC-SM:MANA2015, M-M-N:MMS:2011, M-M:MathNach2010, M-M-N:Springerbook:2013, C-H:2013}.     

The formal asymptotic expansion (\ref{x oustdie1 D_m farmain}) was already given in \cite{Ramm-AG:Springer2005, RAMM:2007}.
The results in Theorem \ref{Maintheorem-ac-small-sing} and in Corollary \ref{Maintheorem-ac-small-sing-1} provide the rigorous justification of the approximation of
the scattered field taking into account all the involved parameters: the maximum diameter $a$, the number $M$, the minimum distance $d$ and the surface impedance $\lambda_m$'s
including the regime defined by $M=M(a):=a^{-s},\; d:=d(a):=a^{t},\; \lambda_m:=\lambda_m(a):=\lambda_{m,0}a^{-\beta}$. We characterized the set of parameters $s, t$ and $\beta$ 
where the approximation makes sense and we provided the approximation with explicit error estimates in terms of these parameters. These approximation formulas can be used for 
different purposes:
\begin{enumerate}
 \item First, they can be used for imaging where the small anomalies are modeled by the small bodies. There is a large literature on the mathematical imaging of small
anomalies, see for instance \cite{A-K:2007,AH-BE-GJ-KH-LH-WA:Book:PSIAM18462015} and the references therein.  Compared to those results, we can allow the small bodies to be dense 
and very close, since their number $M$ can be very large and the minimum distance $d$ can be as small as we want.  This last property means that we deal with the mesoscale regime, where $\frac{a}{\min_{j\neq m}dist(z_j, z_m)}\sim 1$ and $\min_{j\neq m}dist(z_j, z_m)$ is the minimum distance between the centers of the small scatterers, compare with \cite{M-M:MathNach2010,M-M-N:MMS:2011}. 

\item Second, these approximations can be used for
the design of new indices of refraction. When the small scatterers are periodically distributed, this can be justified  by homogenization, see \cite{B-L-P:1978, M-K:2006, J-K-O:1994, DC-FM:Book:BB1997}. Using the approximations in (\ref{x oustdie1 D_m farmain}) and (\ref{x oustdie1 D_m farmain-1}), we do not need such periodicity in the distribution of the small scatterers. This observation has been already made in \cite{RAMM:2007} with quite formal computations. Let us observe that in our approximations 
in Theorem \ref{Maintheorem-ac-small-sing} and in Corollary \ref{Maintheorem-ac-small-sing-1}, we allow the surface impedance to have negative imaginary parts. With this 
freedom of taking the surface impedance, we can generate a large class of 
indices of refraction with possible applications to the acoustic metamaterials.
Details on this issue are reported in \cite{BA-DPC-KMSM:JMAA:2015} and
\cite{C-S:2016-2}.
Let us emphasize that we provide those results with explicit error estimates. 
This gives more credits to the feasibility of the design process.

\item Third, these approximations can be used to extract the values of the index of refraction $n(x), x \in \Omega,$ from the farfields corresponding to few incident directions. The idea of adding small inclusions to deform the medium and then extract the field inside the support of the medium is already described and used in the literature, see for instance \cite{A-C-D-R-T}. 
The extraction of the index of refraction (or other coefficients) from the internal measurements, related also to the hybrid imaging methods, 
has recently attracted the attention of many authors, see for instance 
\cite{Ammari:2008, A-B-G-S:2012, B-B-M-T, B:2014} and the references therein. 
One of the main difficulties is to handle the possible zeros of the
total field. What we propose is to deform the medium using multiple 
(and close) inclusions instead of only single ones. 
In this case, what we derive from the asymptotic expansions are the 
internal values of the Green's function and not only the total fields.
Finally, the values of the index of refraction can be extracted from the 
singularities of these Green's function. Hence, we avoid the problems coming
from the zeros of the internal fields. These arguments are reported in 
\cite{C-S:2016-3}. 

 \end{enumerate}
\bigskip 

To finish this introduction, we discuss the two 'extreme' cases given by the 
Neumann and the Dirichlet boundary conditions. 
In the former case, $\lambda_m$'s are all zero. 
In the results stated above, we see that $\lambda_m$'s are not 
 allowed to be simultaneously zero. Otherwise all the coefficients $C_m,\;
 m=1, ...,M$, will be zero and hence all the 
 constants $Q_m,\;
 m=1, ...,M$, vanish and the expansion (\ref{x oustdie1 D_m farmain}) will make no sense.
The case where $\lambda_m=0,\;
 m=1, ...,M$, is quite tedious but interesting.  
\begin{itemize}
 \item Technically, we see in (\ref{Q_mintdbl-3}), for instance, that in this 
 case we need to go to the higher order in the expansion. Doing that requires quite tedious computations remembering 
that we are taking into account all the parameters describing the scatterers 
($M, a$ and $d$). 
\item This particular case is interesting because in this case the dominant coefficients are defined by {\it{matrices}} 
and not vectors (as the vector $(C_1, C_2, ..., C_M)$ in the case where $\lambda_m \neq 0$). But this is not a surprise since the dominant 
term of the expansion of the far-fields are the Foldy-Lax fields modeling the interaction of the 
multiple point-like scatterers (given here by the 'centers' of the scatterers).
For hard scatterers, we talk about anisotropic 
interaction between the scatterers, see the book \cite{Martin:2006} 
for more information about this issue, contrary to the impedance case where the 
interaction is isotropic. Due to this 'anisotropic' character of the dominant 
term, the equivalent medium (when we distribute a cluster of small scatterers 
with Neumann boundary conditions) is characterized by
a divergence form Helmholtz model where the coefficient appearing in the higher order derivative
is a matrix (defined by the (scaled) matrix-coefficients appearing in the 
dominant term of the expansions).

\end{itemize} 
These arguments need to be mathematically justified and quantified. 
The approximation in the Dirichlet case is discussed in \cite{DPC-SM:MMS2013}
where it is justified under the condition that $\sqrt{M-1}\frac{a}{d}$ is 
bounded, by a constant depending only some a priori bounds,  which means, 
in the scales we use here, that $0\leq s\leq 2-2t$, or 
$t \leq \frac{2-s}{2} $. 
Since the small obstacles are distributed in a bounded domain $\Omega$, then we have the natural
condition on their number $M= O(d^{-3})$, i.e. $\frac{s}{3}\leq t$. 
In the present work, this last condition is the only one we (naturally) impose
on $t$. Hence, compared to \cite{DPC-SM:MMS2013}, we allow here the small scatterers to be as 
close as we want since $t$ can be as large as we want, i.e. 
we cover the mesoscale regime. However, we believe that the conditions used in \cite{DPC-SM:MMS2013}
can be improved to be comparable to the ones we impose here.
\bigskip

 The rest of the paper is devoted to prove 
 Theorem \ref{Maintheorem-ac-small-sing} and 
 Corollary \ref{Maintheorem-ac-small-sing-1}. In section 2, we describe briefly
 the main steps of the proof of Theorem \ref{Maintheorem-ac-small-sing}.
The detailed proof of Theorem \ref{Maintheorem-ac-small-sing} is done 
in  section 3 while the one of Corollary \ref{Maintheorem-ac-small-sing-1} is given in section 4. The justification of Remark \ref{Remark-spherical-scatterers} is discussed in section 5. In the appendix, we deal 
with invertibility of the algebraic system (\ref{fracqcfracmain}) 
(and similarly for (\ref{fracqcfracmain-1})).

\section{A brief description of the proof of Theorem \ref{Maintheorem-ac-small-sing}}

Let us here describe very briefly the main steps of the proof of Theorem \ref{Maintheorem-ac-small-sing}. Firts of all, the scattering problem has a unique solution and it can be 
represented via single layer potentials 
 \begin{equation}\label{SL-rep}
  U^{t}(x)=U^{i}(x)+\sum_{m=1}^{M}\int_{\partial D_m}\Phi_\kappa(x,s)\sigma_{m} (s)ds,~x\in\mathbb{R}^{3}\backslash\left(\mathop{\cup}_{m=1}^M \bar{D}_m\right), 
\end{equation}
where $\sigma:=\left(\sigma_1,\dots,\sigma_M\right)^T$ satisfies the corresponding
system of integral equations. We show that this system of integral equation is invertible
 in the space $\prod\limits_{m=1}^{M}L^2(\partial D_m)$ and that the solution of the scattering
 problem is unique under some natural conditions on the eventual negative imaginary parts of the surface impedance.

We divide the rest of the analysis into few steps:
\begin{enumerate}

\item We have the following a priori estimate of the densities $\sigma_m$'s. 
If $s \leq 2-\beta$ then $\Vert \sigma_m \Vert_{L^2(\partial D_m)} \leq  C a^{1-\beta}$
where $C$ depends only on the Lipschitz characters of 
$B_m,\; m=1, ..., M$. The delicate work here is to derive the precise scaling of the corresponding
boundary integral operators in the appropriate boundary Sobolev spaces 
(taking into account the three parameters $s, t$ and $\beta$).

\item From the representation $
 U^{s}(x)=\sum_{m=1}^{M}\int_{\partial D_m}\Phi_\kappa(x,s)\sigma_{m} (s)ds,\text{ for }x\in\mathbb{R}^{3}\backslash\left(\mathop{\cup}\limits_{m=1}^M \bar{D}_m\right)
$, we deduce, using the above a priori estimates on $\sigma_m$'s, that
\begin{eqnarray}\label{Brief-xfarawayimpnt}
U^{\infty}(\hat{x})&=&\sum_{m=1}^{M}\int_{\partial D_m}e^{-i\kappa\hat{x}\cdot s}\sigma_{m} (s)ds\nonumber\\
 &=&\sum_{m=1}^{M}\left(\int_{\partial D_m}e^{-i\kappa\hat{x}\cdot\,z_{m}}\sigma_{m}(s)ds+\int_{\partial D_m}[e^{-i\kappa\hat{x}\cdot\,s}-e^{-i\kappa\hat{x}\cdot\,z_{m}}]\sigma_{m}(s)ds\right)\nonumber\\
&=&\sum_{m=1}^{M}e^{-i\kappa\hat{x}\cdot\,z_{m}}\tilde{Q}_m+O(\kappa\, M a^{3-\beta})
\end{eqnarray}
where $\tilde{Q}_m:=\int_{\partial D_m}\sigma_m(s)ds$.

\item To estimate the terms $\tilde{Q}_m$, we use the boundary conditions. 
For $s_m\in \partial D_m$, using the impedance boundary condition \eqref{acgoverningsupport}, we have
\begin{eqnarray}\label{Breif-Q_mintdbl}
 0&=&\frac{\partial U^{t}}{\partial \nu_m}(s_m)+\lambda_m U^{t} (s_m)=
-\frac{\sigma_{m} (s_m)}{2}+\int_{\partial D_m}\frac{\partial\Phi_\kappa}{\partial \nu_m}(s_m,s)\sigma_{m} (s)ds+\sum_{\substack{j=1\\j\neq m}}^{M}\int_{\partial D_j}\frac{\partial\Phi_\kappa}{\partial \nu_m}(s_m,s)\sigma_{j} (s)ds\nonumber\\
&& +\lambda_m\sum_{j=1}^{M}\int_{\partial D_j}\Phi_\kappa(s_m,s)\sigma_{m} (s)ds+\frac{\partial U^i}{\partial \nu_m}{(s_m)}+\lambda_m U^i{(s_m)}\nonumber\\
\end{eqnarray}
Integrating the above on $\partial D_m$, we obtain
\begin{equation}%\label{Q_mintdbl-1}
\begin{split}
-\frac{1}{2}\int_{\partial D_m}\sigma_{m} (s_m) d{s_m}+\int_{\partial D_m}\left(\int_{\partial D_m}\frac{\partial\Phi_\kappa}{\partial \nu_m}(s_m,s) d{s_m}\right) \sigma_{m} (s)ds &+\sum_{\substack{j=1\\j\neq m}}^{M}\int_{\partial D_j}\left(\int_{\partial D_m}\frac{\partial\Phi_\kappa}{\partial \nu_m}(s_m,s)d{s_m}\right) \sigma_{j} (s)ds
 \nonumber\\
 +\lambda_m\sum_{j=1}^{M}\int_{\partial D_j}\left(\int_{\partial D_m}\Phi_\kappa(s_m,s) d{s_m}\right) \sigma_{m} (s)ds&=-\int_{\partial D_m}\frac{\partial U^i}{\partial \nu_m}{(s_m)} d{s_m}-\int_{\partial D_m}\lambda_m U^i{(s_m)}d{s_m}\nonumber\\
\end{split}
\end{equation}
It can be rewritten as
\begin{eqnarray}\label{Brief-Q_mintdbl-1}
\begin{split}
-\frac{1}{2}\tilde{Q}_m&+\underbrace{\int_{\partial D_m}\left(\int_{\partial D_m}\frac{\partial\Phi_{0}}{\partial \nu_m}(s_m,s) d{s_m}\right) \sigma_{m} (s)ds}_{=: A} +
\underbrace{\sum_{\substack{j\neq m}}^{M}\int_{\partial D_j}\left(\int_{\partial D_m}\frac{\partial\Phi_\kappa}{\partial \nu_m}(s_m,s)d{s_m}\right) \sigma_{j} (s)ds}_{=: B}
 \nonumber\\
 &+\lambda_m\underbrace{\int_{\partial D_m}\left(\int_{\partial D_m}\Phi_\kappa(s_m,s) d{s_m}\right) \sigma_{m} (s)ds}_{=: C}+\lambda_m\underbrace{\sum_{j\neq m}^{M}\int_{\partial D_j}\left(\int_{\partial D_m}\Phi_\kappa(s_m,z_j) d{s_m}\right) \sigma_{j} (s)ds}_{=: D}
 \nonumber\\
 &=-\int_{\partial D_m}\frac{\partial U^i}{\partial \nu_m}{(s_m)} d{s_m}-\int_{\partial D_m}\lambda_m U^i{(s_m)}d{s_m}+A^{\prime}+\lambda_m D^{\prime},
\end{split}\\
\end{eqnarray}
with
\begin{eqnarray}
A{^\prime}&:=&\int_{\partial D_m}\left(\int_{\partial D_m}\left[\frac{\partial\Phi_\kappa}{\partial \nu_m}(s_m,s)-\frac{\partial\Phi_0}{\partial \nu_m}(s_m,s)\right] d{s_m}\right) \sigma_{m} (s)ds \label{Brief-defofAprm}\\
D{^\prime}&:=&\sum_{j\neq m}^{M}\int_{\partial D_j}\left(\int_{\partial D_m}\left[\Phi_\kappa(s_m,s)-\Phi_\kappa(s_m,z_j)\right] d{s_m}\right) \sigma_{j} (s)ds \label{Brief-defofDprm}. 
\end{eqnarray}
Using the a priori estimate of $\sigma_m$'s, the singularities of the fundamental 
solutions $\Phi_\kappa$ and the harmonicity of $\Phi_0$, we derive the estimates:  
\begin{eqnarray}\label{Brief-approximating-A}
A=-\frac{1}{2}\tilde{Q}_m,\; \mbox{   } A{^\prime}=O\left(\kappa^2a^{4-\beta}\right),
\; B= O\left(2\kappa^2\frac{a^{5-\beta}}{d^{2\alpha}}[\frac{6}{d^{\alpha}}+
7]\right),\;  \mbox{  } 
C=O\left( a^{3-\beta} \right),
  \end{eqnarray}
 \begin{eqnarray}\label{Brief-approximating-Dpr}
D=\sum_{\substack{j\neq m}}^{M}\Phi_\kappa(z_m,z_j) \tilde{Q}_j |\partial D_m|+
O\left(2\frac{a^{5-\beta}}{d^{2\alpha}}\left[7\kappa+\frac{6\kappa+13}{d^{\alpha}} \right] \right)
\end{eqnarray} 
   and
\begin{eqnarray}\label{Brief-approximating-Dpr-1} 
D{^\prime}=O\left(2\frac{a^{5-\beta}}{d^{2\alpha}}\left[7\kappa+\frac{6\kappa+13}{d^{\alpha}} 
\right]\right).
\end{eqnarray}
Here the parameter $\alpha$ is introduced to count the number of small scatterers
surrounding a given and fixed one, see Fig \ref{fig:1-acsmall} and the discussion before. From (\ref{Brief-Q_mintdbl-1}), we obtain the approximation below

\begin{eqnarray}\label{Brief-Q_mintdbl-3}
-\tilde{Q}_m&+&\sum_{\substack{j\neq m}}^{M}\Phi_\kappa(z_m,z_j)\lambda_m|\partial D_m|\tilde{Q}_j\nonumber\\
&=&-\lambda_m |\partial D_m| e^{i\kappa\theta\cdot{z_m}}+O\left((\vert\lambda_m\vert+\kappa)\kappa a^3\right)+\lambda_m O\left(a^{3-\beta}\right),\nonumber\\
&&\quad+O\left(\kappa^2 a^{4-\beta}\right)+O\left(2\kappa^2\frac{a^{5-\beta}}{d^{2\alpha}}[\frac{6}{d^{\alpha}}+7]\right)+\lambda_m O\left(2\frac{a^{5-\beta}}{d^{2\alpha}}\left[7\kappa+\frac{6\kappa+13}{d^{\alpha}} \right]\right).% \nonumber\\
\end{eqnarray}

\par 
Dividing by $C_m:=-\lambda_m|\partial D_m|$ and since 
$\lambda_m=O(a^{-\beta})$, and then $C_m=O(a^{2-\beta})$, we can rewrite the above system as
\begin{eqnarray}\label{Brief-Q_mintdbl-4f}
\frac{\tilde{Q}_m}{C_m}
&=&- e^{i\kappa\theta\cdot{z_m}}-\sum_{\substack{j\neq m}}^{M}
C_j \Phi_\kappa(z_m,z_j)\frac{\tilde{Q}_j}{C_j}+ O
\left( a^{1-\beta}+\frac{a^{3-\beta}}{d^{3\alpha}}\right).
\end{eqnarray}

\item Let now the vector $Y:=(Y_1, Y_2, ..., Y_M)$ be the solution of the Foldy-Lax algebraic 
system
\begin{equation}\label{Brief-FL-Y}
 Y_m =-e^{i\kappa\theta\cdot{z_m}}-\sum_{\substack{j=1 \\ j\neq m}}^{M}C_j 
 \Phi_\kappa(z_m,z_j)Y_j,\; \mbox{ } m=1, ...,\; M.
\end{equation}
We show that the algebraic system (\ref{Brief-FL-Y}) is invertible under general 
condition on $s$ and $\beta$ and derive an error estimate. Based on this error estimate,
we deduce from (\ref{Brief-Q_mintdbl-4f}) and (\ref{Brief-FL-Y}) that

\begin{eqnarray}\label{Brief-unncmaybe}
\sum_{m=1}^{M}\vert {\tilde{Q}_m}-C_mY_m\vert &=& O\left(M{a^{2-\beta}}
\left( a^{1-\beta}+\frac{a^{3-\beta}}{d^{3\alpha}}\right)\right).
\end{eqnarray}

\item The proof ends by plugging (\ref{Brief-unncmaybe}) and 
(\ref{Brief-FL-Y}) in
(\ref{Brief-xfarawayimpnt}) and setting $Q_m:=C_m Y_m,\; \mbox{ } m=1, ...,\; M$.
\end{enumerate}
\bigskip
% % % % % \%\%\%\%\%\%\%\%\%\%\%\%\%\%\%\%\%\%\%\%\%\%\%\%\%\%\%\%\%\%\%\%\%\%\%\%\%\%\%\%\%\%\%\%\%\%\%\%\%
\section{The detailed proof of Theorem \ref{Maintheorem-ac-small-sing}}\label{Proof of Theorem Small}

\subsection{The representation via layer potential}\label{SLPR-1}
We start with the following proposition on the solution of the problem (\ref{acimpoenetrable}-\ref{radiationc}) via the layer potential representation.
% % % \ ~ \ 
% % % \\
\begin{proposition}\label{existence-of-sigmas}%\textit{(Existence of $\sigma_m$'s)}
Assume that the negative part of the imaginary part of $\lambda_m$ are small enough. In addition, suppose that $\kappa^2$ is not a Dirichlet 
\footnote{This last condition is satisfied for every $\kappa$ such that $\kappa\leq\kappa_{\max}$ and 
$a<\frac{1}{\kappa_{\max}}\sqrt[3]{\frac{4\pi}{3}}{\rm j}_{1/2,1}$.
Here ${\rm j}_{1/2,1}$ is the 1st positive zero of the Bessel function ${\rm J}_{1/2}$.} eigenvalue of the Laplacian in $D_m$, for 
$m=1, ...,M$. Then \footnote{ The result of this proposition is valid regardless of the smallness of the obstacles $D_m$'s nor the conditions on $M$, $d$ and $\lambda_m$'s.
The emphasize is on the possibility to deal with surface impedance eventually having negative imaginary parts.}
for $m=1,2,\dots,M$, there exists $\sigma_m\in L^2(\partial D_m)$ such that the problem (\ref{acimpoenetrable}-\ref{radiationc}) has one and a unique solution  and it is 
of the form
 \begin{equation}\label{qcimprequiredfrm1}
  U^{t}(x)=U^{i}(x)+\sum_{m=1}^{M}\int_{\partial D_m}\Phi_\kappa(x,s)\sigma_{m} (s)ds,~x\in\mathbb{R}^{3}\backslash\left(\mathop{\cup}_{m=1}^M \bar{D}_m\right), 
\end{equation}
% with $\Phi_\kappa(x,y):=\frac{e^{i\kappa|x-y|}}{4\pi|x-y|},$ for all $x,y$ in $\mathbb{R}^3$. This solution is unique.
\end{proposition}
\begin{proofn}{\it{of Proposition \ref{existence-of-sigmas}}.}
  We look for the solution of the problem (\ref{acimpoenetrable}-\ref{radiationc}) of the form \eqref{qcimprequiredfrm1}, 
  then from the impedance boundary condition \eqref{acgoverningsupport}, we obtain
\begin{eqnarray}\label{qcimprequiredfrm1bd}
-\frac{\sigma_{j} (s_j)}{2}+\int_{\partial D_j}\frac{\partial\Phi_\kappa(s_j,s)}{\partial \nu_j(s_j)}\sigma_{j} (s)ds+\sum_{\substack{m=1\\m\neq j}}^{M}\int_{\partial D_m}\frac{\partial\Phi_\kappa(s_j,s)}{\partial \nu_j(s_j)}\sigma_{m} (s)ds&&\\
+\lambda_j\sum_{m=1}^{M}\int_{\partial D_m}\Phi_\kappa(s_j,s)\sigma_{m} (s)ds&=&-\frac{\partial U^i{(s_j)}}{\partial \nu_j(s_j)}-\lambda U^i{(s_j)},\,\forall s_j\in \partial D_j,\, j=1,\dots,M.\nonumber
\end{eqnarray}
One can write it in a compact form as 
$(-\frac{1}{2}\textbf{I}+DL^{*}+DK^{*}+\lambda(L+K))\sigma=-(\partial_\nu+\lambda) U^{In}$ with $\partial_\nu:=({\partial_\nu}_{mj})_{m,j=1}^{M}$, $\lambda:=(\lambda_{mj})_{m,j=1}^{M}$, $DL^{*}:=(DL^{*}_{mj})_{m,j=1}^{M}$, $DK^{*}:=(DK^{*}_{mj})_{m,j=1}^{M}$,
$L:=(L_{mj})_{m,j=1}^{M}$ and $K:=(K_{mj})_{m,j=1}^{M}$, where

\begin{eqnarray}\label{definition-L_K}
\hspace{-.3cm}\textbf{I}_{mj}=\left\{\begin{array}{ccc}
            I,\, \text{Identity operator} & m=j\\
            0,\, \text{zero~operator}~~~~~~ & else
           \end{array}\right.,\,
&
DL^{*}_{mj}=\left\{\begin{array}{ccc}
            \mathcal{D}^{*}_{mj} & m=j\\
            0 & else
           \end{array}\right.,
&  DK^{*}_{mj}=\left\{\begin{array}{ccc}
            \mathcal{D}^{*}_{mj} & m\neq j\\
            0 & else
           \end{array}\right., 
\end{eqnarray}
\begin{eqnarray}\label{definition-L_K-1}
{\partial_\nu}_{mj}=\left\{\begin{array}{ccc}
            \partial_{\nu_{m}} & m=j\\
            0 & else
           \end{array}\right.,
\quad
\lambda_{mj}=\left\{\begin{array}{ccc}
            \lambda_{m} & m=j\\
            0 & else
           \end{array}\right.,
&&\,
L_{mj}=\left\{\begin{array}{ccc}
            \mathcal{S}_{mj} & m=j\\
            0 & else
           \end{array}\right.,
 \quad K_{mj}=\left\{\begin{array}{ccc}
            \mathcal{S}_{mj} & m\neq j\\
            0 & else
           \end{array}\right., \,
\end{eqnarray}

$U^I=U^I(s_1,\dots,s_M):=\left(U^i(s_1),\dots,U^i(s_M)\right)^T$ and $\sigma=\sigma(s_1,\dots,s_M):=\left(\sigma_1(s_1),\dots,\sigma_M(s_M)\right)^T$. 
Here, for the indices  $m$ and $j$ fixed, $\mathcal{S}_{mj}$ is the integral operator acting as

\begin{eqnarray}\label{defofSmjed}
 \mathcal{S}_{mj}(\sigma_j)(t):=\int_{\partial D_j}\Phi_\kappa(t,s)\sigma_j(s)ds,\quad t\in\partial D_m,
\end{eqnarray}
and $\mathcal{D}^{*}_{mj}$ is the adjoint of the integral operator  defined by,

\begin{eqnarray}\label{defofDmjed}
 \mathcal{D}_{mj}(\sigma_j)(t):=\int_{\partial D_j}\frac{\partial\Phi_\kappa(t,s)}{\partial \nu_m(s)}\sigma_j(s)ds,\quad t\in\partial D_m. 
\end{eqnarray}

Then the operator $\mathcal{D}^{*}_{mm}:L^2(\partial D_m)\rightarrow L^2(\partial D_m)$ is adjoint of the double layer operator $\mathcal{D}_{mm}:L^2(\partial D_m)\rightarrow L^2(\partial D_m)$, defined by,
\begin{eqnarray}\label{defofDmm}
 \mathcal{D}_{mm}(\sigma_m)(t):=\int_{\partial D_m}\frac{\partial\Phi_\kappa(t,s)}{\partial \nu_m(s)}\sigma_j(s)ds,\quad t\in\partial D_m,
\end{eqnarray}
 and the operator $-\frac{1}{2}I+\mathcal{D}^{*}_{mm}:L^2(\partial D_m)\rightarrow L^2(\partial D_m)$ is isomorphic and hence
 Fredholm with zero index. For $m\neq j$, $\mathcal{D}^{*}_{mj}:L^2(\partial D_j)\rightarrow L^2(\partial D_m)$ 
is compact, %when $\partial D_m$ has a Lipschitz regularity, 
see \cite[Theorem 4.1]{MD:InteqnsaOpethe1997}.\footnote{Observe that $-\frac{1}{2}I+\mathcal{D}^{*}_{mm}$ is the adjoint of $-\frac{1}{2}I+\mathcal{D}_{mm}$. 
Hence $-\frac{1}{2}I+\mathcal{D}^{*}_{mm}$ is Fredholm if we show that $-\frac{1}{2}I+\mathcal{D}_{mm}$ is. In \cite{MD:InteqnsaOpethe1997}, this last property is proved 
for the case $\kappa=0$ and by a perturbation argument, we can obtain the same results for every $\kappa$. Using the condition on $\kappa^2$, we deduce that 
$-\frac{1}{2}I+\mathcal{D}^{*}_{mm}$ is an isomorphism.}

%It is obtained using the Faber Krahn inequality mentioned in the paper 'eigenvalues of the Laplacian acting on functions of mean zero with constant boundary values' by Antonio Greco. 
%The calculation is $\kappa^2<\lambda_1(\Omega^*)=\left(\frac{\frac{4}{3}\pi}{a^3}\right)^{\frac{2}{3}}{\rm j}^2_{1/2,1}$
 
Also notice that  $\mathcal{S}_{mj}:L^2(\partial D_j)\rightarrow L^2(\partial D_m)$ 
is compact.
% Remark here that, for the scattering by single obstacle $DK^{*}$ and $K$  are zero  operators. %In our case  we consider $r=0$. 
So, $(-\frac{1}{2}\textbf{I}+DL^{*}+DK^{*}+\lambda(L+K)):\prod\limits_{m=1}^{M}L^2(\partial D_m)\rightarrow \prod\limits_{m=1}^{M}L^2(\partial D_m)$ is Fredholm with zero index.  
We induce the product of spaces by the maximum of the norms of the space.
To show that $(-\frac{1}{2}\textbf{I}+DL^{*}+DK^{*}+\lambda(L+K))$ is invertible it is enough to show that it is injective. i.e. $(-\frac{1}{2}\textbf{I}+DL^{*}+DK^{*}+\lambda(L+K))\sigma=0$ implies $\sigma=0$.
We write 
$$\tilde{U}(x)=\sum_{m=1}^{M}\int_{\partial D_m}\Phi_\kappa(x,s)\sigma_m(s)ds, \mbox{ in } \mathbb{R}^3\backslash 
\left(\mathop{\cup}_{m=1}^{M}\bar{D}_m\right)$$
and 
$$\tilde{\tilde{U}}(x)=\sum_{m=1}^{M}\int_{\partial D_m}\Phi_\kappa(x,s)\sigma_m(s)ds, \mbox{ in } \mathop{\cup}_{m=1}^{M}D_m.$$

Then $\tilde{U}$ satisfies $\Delta\tilde{U}+\kappa^2\tilde{U}=0$ for $x\in\mathbb{R}^{3}\backslash\left(\mathop{\cup}\limits_{m=1}^{M}\bar{D}_m\right)$,  
with S.R.C and $(\partial_{\nu_m}+\lambda_m)\tilde{U}(x)=0$ on $\mathop{\cup}\limits_{m=1}^{M}\partial D_m$.  

Let us assume for the moment that the exterior problem has a unique solution, then we deduce that $\tilde{U}=0$ in $\mathbb{R}^{3}\backslash\left(\mathop{\cup}\limits_{m=1}^{M}\bar{D}_m\right)$.
Due to the continuity of the single layer potentials, we deduce that $\tilde{\tilde{U}}=0$ on $\mathop{\cup}\limits_{m=1}^{M}\partial D_m$. 
In addition, we know that $\Delta\tilde{\tilde{U}}+\kappa^2\tilde{\tilde{U}}=0$ for $x\in\mathop{\cup}\limits_{m=1}^{M}D_m$. 
From the condition on $\kappa^2$, we deduce that $\tilde{\tilde{U}}=0$ in $\left(\mathop{\cup}\limits_{m=1}^{M}D_m\right)$. 

By the jump relations, we have
\begin{eqnarray}\label{nbc1}
\frac{\partial\tilde{U}}{\partial \nu}(x)+\lambda_m \tilde{U} (x)
=0&\Longrightarrow&(\mathbf{K}^{*}\sigma_m)(x)-\frac{\sigma_m(x)}{2}+\sum\limits^{M}_{\substack{j=1\\j\neq m}}(\mathcal{D}^{*}_{mj}+\lambda_m \mathcal{S}_{mj})(\sigma_j)(x)=0
\end{eqnarray} 
and 
\begin{eqnarray}\label{nbc2}
\frac{\partial\tilde{\tilde{U}}}{\partial \nu}(x)+\lambda_m \tilde{\tilde{U}} (x)=0&\Longrightarrow&(\mathbf{K}^{*}\sigma_m)(x)+\frac{\sigma_m(x)}{2}+\sum\limits^{M}_{\substack{j=1\\j\neq m}}(\mathcal{D}^{*}_{mj}+\lambda_m \mathcal{S}_{mj})(\sigma_j)(x)=0
\end{eqnarray}
 for $x\in \partial D_m$ and for $m=1,\dots,M$. Here, $\mathbf{K}^{*}$ is the adjoint of the double layer operator $\mathbf{K}$,
\begin{eqnarray}\label{nbc23} 
(\mathbf{K}\sigma_m)(x):=\int_{\partial D_m}\frac{\partial}{\partial \nu_s}\Phi_\kappa(x,s) \sigma_m(s) ds,\mbox{ for } m=1,\dots,M. 
\end{eqnarray}
Difference between \eqref{nbc1} and \eqref{nbc2} provides us, $\sigma_m=0$ for all $m$. \\
 %By the jump relations, we see that $\sigma_m=0$ for $m=1,\dots,M$. Hence, the operator $\mathcal{S}:=L+K$ is invertible. 
We conclude then that $-\frac{1}{2}\textbf{I}+DL^{*}+DK^{*}+\lambda(L+K)=:-\frac{1}{2}\textbf{I}+\mathcal{D}^{*}+\lambda \mathcal{S}:
\prod\limits_{m=1}^{M}L^2(\partial D_m)\rightarrow \prod\limits_{m=1}^{M}L^2(\partial D_m)$ is invertible.
% and so $\sigma=\mathcal{S}^{-1}U^I\in  \prod\limits_{m=1}^{M}L^{2}(\partial D_m)$ exists.
\bigskip

We need now to show that the exterior problem $\Delta\tilde{U}+\kappa^2\tilde{U}=0$ for $x\in\mathbb{R}^{3}\backslash\left(\mathop{\cup}\limits_{m=1}^{M}\bar{D}_m\right)$,  
with S.R.C and $(\partial_{\nu_m}+\lambda_m)\tilde{U}(x)=0$ on $\mathop{\cup}\limits_{m=1}^{M}\partial D_m$, has a unique solution. This result is known under the condition 
$\Im \lambda_m \geq 0$ on $\mathop{\cup}\limits_{m=1}^{M}\partial D_m$, see 
\cite{C-K:1983} for instance. To relax this positivity condition and consider 
$\Im \lambda_m < 0$ for some $m$'s, we proceed as follows. 
Set $f:=(f_1,...,f_M)$ with $f_m:=i(\Im \lambda)_- \tilde{U}$ on $\partial D_m$ where 
$(\Im \lambda)_-:=\max_{\mathop{\cup}\limits_{m=1}^{M}\partial D_m}\{-\Im \lambda_m,0\}$. Hence $\tilde{U}$ satisfies 
$\Delta\tilde{U}+\kappa^2\tilde{U}=0$ for $x\in\mathbb{R}^{3}\backslash\left(\mathop{\cup}\limits_{m=1}^{M}\bar{D}_m\right)$,  
with S.R.C and $(\partial_{\nu_m}+\lambda_m+i(\Im\lambda)_-)\tilde{U}(x)=f_m$ on $\mathop{\cup}\limits_{m=1}^{M}\partial D_m$. Since now $\Im(\lambda_m+i(\Im \lambda)_-)>0$, then 
this last problem has a unique solution and it can be represented via single layer potentials $\sum^M_{m=1} S(\psi_j)$. Taking the normal trace on $\partial D_m$'s, 
we deduce that $\left(-\frac{1}{2}\textbf{I}+DL^{*}+DK^{*}+(\lambda+i(\Im\lambda)_-\mathbf{I})(L+K)\right)\psi=f$, where $\psi:=(\psi_1,...\psi_M)$. But from the form of $f$, we obviously have
$f=i(\Im\lambda)_-(L+K)\psi$. Hence 
\begin{equation}\label{spectral problem-boundary}
(-\frac{1}{2}\textbf{I}+DL^{*}+DK^{*}+\lambda(L+K))\psi=0.
\end{equation}
Let us first consider, for simplicity, only one scatterer and assume that the surface impedance is a constant 
 $\lambda:=\lambda^r + i\;\lambda^i$. Then (\ref{spectral problem-boundary}) reduces to 
\begin{equation}\label{boundary-equation}
 (-\frac{1}{2}I+\mathbf{K}^{*}+\lambda S)\psi=0.
\end{equation}

1. If $\lambda^i \geq 0$, then, as usual in the scattering theory \cite{C-K:1983}, we derive 
that $S\psi=0$ 
and then $\psi=0$.
\bigskip

2. How about $\lambda^i<0?$ Since the operator $-\frac{1}{2}I+\mathbf{K}^* 
+ \lambda^r S$ is 
invertible in the $L^2(\partial D_m)$ spaces, then the equation (\ref{boundary-equation}) can be reduced to
\begin{equation}\label{boundary-equation-2}
 (-\frac{1}{2}I+\mathbf{K}^*+ \lambda^r S)^{-1}S\psi=-(i\lambda^i)^{-1}\psi=i
 (\lambda^i)^{-1}\psi
\end{equation}
i.e. $(\psi, i(\lambda^i)^{-1})$ is an eigenelement of the operator 
$(-\frac{1}{2}I+\mathbf{K}^*+\lambda^r S)^{-1}S$. 
But $(-\frac{1}{2}I+\mathbf{K}^*+\lambda^r S)^{-1}S$ is compact hence it has 
only a discrete set of eigenvalues. 
%\footnote{If the scatterer is a ball fo radius $a$ centered at the origin, indeed for $\kappa$ and $a$ fixed, and for 
%$\lambda:=\frac{1}{a}-i \kappa$, $i\kappa^{-1}$ is an eigenvalue of the 
%related operator $(-\frac{1}{2}I+\mathbf{K}^*+\frac{1}{a} S)^{-1}S$ where 
%$\psi:=S^{-1}\frac{e^{i \kappa r}}{r}$ is a corresponding eigenfunction.}. 
Then if we take the surface impedance $\lambda$ such that 
$i(\lambda^i)^{-1}$ 
is different from these discrete values, then $\psi=0$ and hence $U=0$.
\bigskip

We conclude that the scattering by an obstacle with an impedance type boundary condition modeled by 
a constant $\lambda^r +i \lambda^i$ 
 is well posed as soon as the imaginary part $\lambda^i$ is such that $i(\lambda^i)^{-1}$ is not an eigenvalue of 
the corresponding compact operator $(-\frac{1}{2}I+\mathbf{K}^* +\lambda^r S)^{-1}S$.
\bigskip

In particular, the operator in ({\ref{boundary-equation}}), i.e. $-\frac{1}{2}I+\mathbf{K}^*+\lambda S$, 
can be inverted using the Neumann series if $\lambda$ satisfies 
$
 \vert \lambda^i \vert \Vert (-\frac{1}{2}I+\mathbf{K}^*+ \lambda^r S)^{-1}\Vert \Vert S \Vert<1. 
$ This is of course a stronger condition on $\lambda$ but it is enough for our purposes. 
In addition this Neumann series argument applies smoothly to the case where we have variable 
surface impedance's and multiple scatterers. Indeed, we know that the operator 
$-\frac{1}{2}\textbf{I}+DL^{*}+DK^{*}+(\lambda+i(\Im \lambda)_-\mathbf{I})(L+K):$ 
$\prod\limits_{m=1}^{M}L^2(\partial D_m)\rightarrow \prod\limits_{m=1}^{M}L^2(\partial D_m)$ is invertible. 
Hence if the negative part of the imaginary part of $\lambda$ is small so 
that \footnote{A general condition to inverte 
(\ref{spectral problem-boundary}) is to assume that $-1$ is not an eigenvalue of the compact operator
$(-\frac{1}{2}\textbf{I}+DL^{*}+DK^{*}+(\lambda+i(\Im \lambda)_-\mathbf{I})
(L+K))^{-1}(i(\Im \lambda)_ -(L+K))$.}

\begin{equation}\label{negative-imaginary-lambda}
 (\Im \lambda)_-\Vert L+K \Vert \Vert (-\frac{1}{2}\textbf{I}+DL^{*}+DK^{*}+(\lambda+i(\Im \lambda)_-\mathbf{I})(L+K))^{-1}\Vert <1, 
\end{equation}
then by the Neumann series expansion $-\frac{1}{2}\textbf{I}+DL^{*}+DK^{*}+\lambda(L+K)$ is also invertible and hence $\psi=0$.

\end{proofn}
 %In the similar way, we can show that $L$ is invertible. 
\subsection{An appropriate estimate of the densities $\sigma_m,\,m=1,\dots,M$}\label{SLPR-2}
From the above theorem, we have the following representation of $\sigma$:
\begin{eqnarray}\label{invLplusK} 
 \sigma&=&-(-\frac{1}{2}\textbf{I}+DL^{*}+DK^{*}+\lambda(L+K))^{-1}(\partial_\nu+\lambda) U^{In} \\
       &=&-(-\frac{1}{2}\textbf{I}+DL^{*}+\lambda L)^{-1}(\textbf{I}+(-\frac{1}{2}\textbf{I}+DL^{*}+\lambda L)^{-1}(DK^{*}+\lambda K))^{-1}(\partial_\nu+\lambda) U^{In} \nonumber\\
       &=&-(-\frac{1}{2}\textbf{I}+DL^{*}+\lambda L)^{-1}\sum_{l=0}^{\infty}\left((-\frac{1}{2}\textbf{I}+DL^{*}+\lambda L)^{-1}(DK^{*}+
       \lambda K)\right)^{l}(\partial_\nu+\lambda) U^{In}, \nonumber
\end{eqnarray}
if $ \left\|(-\frac{1}{2}\textbf{I}+DL^{*}+\lambda L)^{-1}(DK^{*}+\lambda K)\right\|<1.$
\noindent
By the assumption $ a<\frac{1}{\kappa_{\max}}\sqrt[3]{\frac{4\pi}{3}}{\rm j}_{1/2,1}$, the operator $-\frac{1}{2}\textbf{I}+DL^{*}$ is invertible.
 We write $-\frac{1}{2}\textbf{I}+DL^{*}+\lambda L=(I+\lambda L(-\frac{1}{2}\textbf{I}+DL^{*})^{-1})(-\frac{1}{2}\textbf{I}+DL^{*})$.
 From the scaling of the operators $L$ and $(-\frac{1}{2}\textbf{I}+DL^{*})$, we show that there exists a constant $\tilde{c}$, depending only on the Lipschitz character of
 the reference bodies $B_m$'s, such that if 
 \begin{equation}\label{cond-lambda}
  \| \lambda \|a =\max_{m}\vert \lambda_{m,0}\vert a^{1-\beta} <\tilde{c},
 \end{equation}
we have $\Vert \lambda L(-\frac{1}{2}\textbf{I}+DL^{*})^{-1} \Vert <1$ and hence the operator $-\frac{1}{2}\textbf{I}+DL^{*}+\lambda L$ is invertible.
The condition (\ref{cond-lambda}) is  verified if $\beta <1$ and $a$ small enough. For the convinience, we denote  $\| \lambda \|$ by $\vert \lambda \vert$.

This implies that
\begin{eqnarray}\label{nrminvLplusK}
 \left\|\sigma\right\| 
                                  &\leq&\frac{\left\|(-\frac{1}{2}\textbf{I}+DL^{*}+\lambda L)^{-1}\right\|}{1-\left\|(-\frac{1}{2}\textbf{I}+DL^{*}+\lambda L)^{-1}\right\|\left\|DK^{*}+\lambda K)\right\|}\left\|(\partial_\nu+\lambda) U^{In}\right\|.
\end{eqnarray}
Here we use the following notations: 
\begin{eqnarray}
\left\|DK^{*}+\lambda K\right\|&:= &\left\|DK^{*}+\lambda K\right\|_{\mathcal{L}\left(\prod\limits_{m=1}^{M}L^{2}(\partial D_m),\prod\limits_{m=1}^{M}L^{2}(\partial D_m)\right)}\nonumber\\
    &\equiv&%\sum\limits_{m=1}^{M}\sum\limits_{j=1}^{M}
 \max\limits_{1\leq m \leq M}\sum_{j=1}^{M}\left\|DK^{*}_{mj}+\lambda_m K_{mj}\right\|_{\mathcal{L}\left(L^{2}(\partial D_j),L^{2}(\partial D_m)\right)}\nonumber\\
    &=&\max\limits_{1\leq m \leq M}\sum_{\substack{j=1\\j\neq\,m}}^{M}\left\|\mathcal{D}^{*}_{mj}+\lambda_m \mathcal{S}_{mj}\right\|_{\mathcal{L}\left(L^{2}(\partial D_j),L^{2}(\partial D_m)\right)},\label{DKSnrm}\\
% % % %     \left\|(L+K)\right\|&:= &\left\|(L+K)\right\|_{\mathcal{L}\left(\prod\limits_{m=1}^{M}L^{2}(\partial D_m),\prod\limits_{m=1}^{M}L^{2}(\partial D_m)\right)}\nonumber\\
% % % %     &\equiv&%\sum\limits_{m=1}^{M}\sum\limits_{j=1}^{M}
% % % %  \max\limits_{1\leq m \leq M}\sum_{j=1}^{M}\left\|(L+K)_{mj}\right\|_{\mathcal{L}\left(L^{2}(\partial D_j),L^{2}(\partial D_m)\right)}\nonumber\\
% % % %     &=&\max\limits_{1\leq m \leq M}\sum_{\substack{j=1}}^{M}\left\|\mathcal{S}_{mj}\right\|_{\mathcal{L}\left(L^{2}(\partial D_j),L^{2}(\partial D_m)\right)},\label{LplusKnrm}\\
\left\|(-\frac{1}{2}\textbf{I}+DL^{*}+\lambda L)^{-1}\right\|
&:=& \left\|(-\frac{1}{2}\textbf{I}+DL^{*}+\lambda L)^{-1}\right\|_{\mathcal{L}\left(\prod\limits_{m=1}^{M}L^{2}(\partial D_m),\prod\limits_{m=1}^{M}L^{2}(\partial D_m)\right)}\nonumber\\
     &\equiv&\max\limits_{1\leq m \leq M}\sum_{j=1}^{M}\left\|{\left((-\frac{1}{2}\textbf{I}+DL^{*}+\lambda_m L)^{-1}\right)_{mj}}\right\|_{\mathcal{L}\left(L^{2}(\partial D_m),L^{2}(\partial D_j)\right)}\nonumber\\
     &=&
\max\limits_{1\leq m \leq M}\left\|(-\frac{1}{2}I+\mathcal{D}^{*}_{mm}+\lambda_m \mathcal{S}_{mm})^{-1}\right\|_{\mathcal{L}\left(L^{2}(\partial D_m),L^{2}(\partial D_m)\right)},\label{invDLSnrm}\\
\left\|\sigma\right\|&:=& \left\|\sigma\right\|_{\prod\limits_{m=1}^{M}L^{2}(\partial D_m)}
    \,\quad\equiv\,\max\limits_{1\leq m \leq M}\left\|\sigma_{m}\right\|_{L^{2}(\partial D_m)},\label{sigmaU^Inrm1}\\ 
    \left\|U^{In}\right\|&:=& \left\|U^{In}\right\|_{\prod\limits_{m=1}^{M}L^{2}(\partial D_m)}
    \,\quad\equiv\,\max\limits_{1\leq m \leq M}\left\|U^{i}\right\|_{L^{2}(\partial D_m)}\label{sigmaU^Inrmder}\\
\mbox{and}\qquad\left\|\partial_\nu U^{In}\right\|&:= &\left\|\partial_\nu U^{In}\right\|_{\prod\limits_{m=1}^{M}L^{2}(\partial D_m)}
    \equiv\,\max\limits_{1\leq m \leq M}\left\|\partial_{\nu_m}U^{i}\right\|_{L^{2}(\partial D_m)}.\label{sigmaU^Inrm}
\end{eqnarray}

In the following proposition, we provide conditions under which $\left\|(-\frac{1}{2}\textbf{I}+DL^{*}+\lambda L)^{-1}(DK^{*}+\lambda K)\right\|<1$ and then estimate $\left\|\sigma\right\|$ via \eqref{nrminvLplusK}.
\begin{proposition}\label{normofsigmastmt}%\textit{(Estimate of $\sigma_m$'s)}
There exists $a_0$ depending only on the 
set of the a priori bounds such that if $a\leq a_0$ and $s\leq 2-\beta$, then we have the following estimate
\begin{eqnarray*} 
 \left\|\sigma_m\right\|_{L^{2}(\partial D_m)}\,\leq\, c \epsilon^{1-\beta}
\end{eqnarray*}
where $c$ is a positive constant depending only on the 
set of the a priori bounds.
 \end{proposition}
\noindent

\textit{Proof of Proposition \ref{normofsigmastmt}}.
 
 Suppose $0<\epsilon\leq1$ and $D_\epsilon:=\epsilon B+z\subset\mathbb{R}^n$. For any functions $f,g$ defined on $\partial D_\epsilon$ and $\partial B$ respectively, we use the notations; %we define 
%$$(f)^\wedge(\xi):=\hat{f}(\xi):=f(\epsilon\xi+z)\mbox{ and } (g)^\vee(x):= \check{g}(x):=g\left(\frac{x-z}{\epsilon}\right).$$ 
\begin{eqnarray}\label{acsmalla-wedgevee-defntn}
 (f)^\wedge(\xi)\,:=\,\hat{f}(\xi)\,:=\,f(\epsilon\xi+z)&\mbox{ and }& (g)^\vee(x)\,:=\, \check{g}(x)\,:=\,g\left(\frac{x-z}{\epsilon}\right).
\end{eqnarray}

 Then for each  $\psi \in L^{2}(\partial D_\epsilon)$, we have
 \begin{equation}\label{habib2*}
 \|\psi \|_{L^{2}(\partial D_\epsilon)}=\epsilon^\frac{n-1}{2} \|\hat{\psi} \|_{L^{2}(\partial B)}
\end{equation}

% \begin{equation}\label{habib1*}
% \epsilon^{\frac{n-1}{2}}\vert\vert \hat{\phi} \vert\vert_{H^1(\partial B)}~\leq ~\vert\vert \phi \vert\vert_{H^{1}(\partial D_\epsilon)}
% ~\leq~\epsilon^{\frac{n-3}{2}}\vert\vert \hat{\phi} \vert\vert_{H^{1}(\partial B)}.
% \end{equation}
\noindent
We divide the rest of the proof of Proposition \ref{normofsigmastmt} into two steps. In the first step, we assume we have a single obstacle and then in the second step we deal with the multiple obstacles case.
% \paragraph{The case of a single obstacle}\label{case of a single obstacle-smallac-sdlp} 
\subsubsection{The case of a single obstacle}\label{case of a single obstacle-smallac-sdlp} 
%In our problem we deal only with the 3d case. i.e. n=3.  
Let us consider a single obstacle $D_\epsilon:=\epsilon B+z$ with unit outword normal $\nu$ to its boundary. %, which has the properties as of $D_m$.  
Then define the operator 
$\mathcal{D}_{ D_\epsilon}:L^2(\partial D_{\epsilon})\rightarrow L^2(\partial D_{\epsilon})$  by
\begin{eqnarray}\label{defofSpartialDe}
\left(\mathcal{D}_{ D_\epsilon}\psi\right) (s)=\int_{\partial D_\epsilon}\frac{\Phi_\kappa(s,t)}{\partial\nu(t)}\psi(t)dt.
\end{eqnarray}
Following the arguments in the proof of Proposition \ref{existence-of-sigmas}, the integral operator 
$-\frac{1}{2}I+\mathcal{D}^{*}_{ D_\epsilon}:L^2(\partial D_{\epsilon})\rightarrow L^2(\partial D_{\epsilon})$ is invertible. 
If we consider the problem (\ref{acimpoenetrable}-\ref{radiationc}) 
in $\mathbb{R}^{3}\backslash \bar{D}_\epsilon$, we obtain
  $$\sigma=(-\frac{1}{2}I+\mathcal{D}^{*}_{ D_\epsilon}+\lambda\mathcal{S}_{ D_\epsilon})^{-1}(\partial_\nu+\lambda) U^{i},$$%\,\text{where}\,DL^{*}+DK^{*}=:\mathcal{D}^{*}_{ D_\epsilon}$$%,~ \mathcal{S}:=\mathcal{S}_{ D_\epsilon}$$ 
and then
\begin{eqnarray}\label{estsigm1}
 \|\sigma\|_{L^2(\partial D_{\epsilon})}\leq \|(-\frac{1}{2}I+\mathcal{D}^{*}_{ D_\epsilon}+\lambda\mathcal{S}_{ D_\epsilon})^{-1}\|_{\mathcal{L}\left(L^2(\partial D_{\epsilon}),L^2(\partial D_{\epsilon})\right)}\|(\partial_\nu+\lambda) U^i\|_{L^2(\partial D_{\epsilon})}.
\end{eqnarray}
We have the following lemma, see \cite[Lemma 2.4 and Lemma 2.15]{DPC-SM:MMS2013}.

\begin{lemma}\label{rep1singulayer}
 Let $\phi, \psi \in L^{2}(\partial D_\epsilon)$. Then, 
 \begin{equation}\label{rep1singulayer1d}
 \mathcal{S}_{ D_\epsilon}\psi=\epsilon~ (\mathcal{S}^\epsilon_B \hat{\psi})^\vee,
\end{equation}

\begin{equation}\label{nrm1singulayer2dd}
 \left\|\mathcal{S}_{ D_\epsilon}\right\|_{\mathcal{L}\left(L^2(\partial D_\epsilon), L^2(\partial D_\epsilon) \right)}=
\epsilon\left\|{\mathcal{S}^\epsilon_B}\right\|_{\mathcal{L}\left(L^2(\partial B), L^2(\partial B) \right)},
\end{equation}

\begin{equation}\label{rep1singulayer1}
 \mathcal{D}^{*}_{ D_\epsilon}\psi= ({\mathcal{D}^{\epsilon^{*}}_B} \hat{\psi})^\vee,
\end{equation}
\begin{equation}\label{rep1sgllayer1}
 \left(-\frac{1}{2}I+\mathcal{D}^{*}_{ D_\epsilon}\right)\psi=  \left(\left(-\frac{1}{2}I+{\mathcal{D}^{\epsilon^{*}}_B} \right)\hat{\psi}\right)^\vee,
\end{equation}
 \begin{equation}\label{rep1sgllayer2}
 {\left(-\frac{1}{2}I+\mathcal{D}^{*}_{ D_\epsilon}\right)}^{-1}\phi= \left({\left(-\frac{1}{2}I+{\mathcal{D}^{\epsilon^{*}}_B}\right)}^{-1} \hat{\phi}\right)^\vee,
\end{equation}

\begin{equation}\label{rep1singulayer2}
 \left\|{\left(-\frac{1}{2}I+\mathcal{D}^{*}_{ D_\epsilon}\right)}^{-1}\right\|_{\mathcal{L}\left(L^2(\partial D_\epsilon), L^2(\partial D_\epsilon) \right)}=
\left\|{\left(-\frac{1}{2}I+{\mathcal{D}^{\epsilon^{*}}_B}\right)}^{-1}\right\|_{\mathcal{L}\left(L^2(\partial B), L^2(\partial B) \right)}
\end{equation}

and hence 
\begin{equation}\label{rep1sgldbllayer1}
 \left(-\frac{1}{2}I+\mathcal{D}^{*}_{ D_\epsilon}+\lambda \mathcal{S}_{ D_\epsilon}\right)\psi=  \left(\left(-\frac{1}{2}I+{\mathcal{D}^{\epsilon^{*}}_B}+\lambda\epsilon \mathcal{S}^\epsilon_B \right)\hat{\psi}\right)^\vee,
\end{equation}
 \begin{equation}\label{rep1sgldbllayer2}
 {\left(-\frac{1}{2}I+\mathcal{D}^{*}_{ D_\epsilon}+\lambda \mathcal{S}_{ D_\epsilon}\right)}^{-1}\phi= \left({\left(-\frac{1}{2}I+{\mathcal{D}^{\epsilon^{*}}_B}+\lambda\epsilon \mathcal{S}^\epsilon_B\right)}^{-1} \hat{\phi}\right)^\vee,
\end{equation}

\begin{equation}\label{rep1singuldblayer2}
 \left\|{\left(-\frac{1}{2}I+\mathcal{D}^{*}_{ D_\epsilon}+\lambda \mathcal{S}_{ D_\epsilon}\right)}^{-1}\right\|_{\mathcal{L}\left(L^2(\partial D_\epsilon), L^2(\partial D_\epsilon) \right)}=
\left\|{\left(-\frac{1}{2}I+{\mathcal{D}^{\epsilon^{*}}_B}+\lambda\epsilon \mathcal{S}^\epsilon_B\right)}^{-1}\right\|_{\mathcal{L}\left(L^2(\partial B), L^2(\partial B) \right)},
\end{equation}

with  ${\mathcal{S}^{\epsilon}_B} \hat{\psi}(\xi):=\int_{\partial B}\Phi^\epsilon(\xi,\eta)\hat{\psi}(\eta) d\eta$, ${\mathcal{D}^{\epsilon^{*}}_B} \hat{\psi}(\xi):=\int_{\partial B}\frac{\partial\Phi^\epsilon(\xi,\eta)}{\partial \nu(\xi)}\hat{\psi}(\eta) d\eta$ and
$\Phi^{\epsilon}(\xi,\eta):=\frac{e^{i\kappa\epsilon|\xi-\eta|}}{4\pi|\xi-\eta|}$.% and $\left(\frac{1}{2}I+\mathcal{D}\right)_{ D_\epsilon~(B)}=\frac{1}{2}I+\mathcal{D}_{ D_\epsilon~(B)}$.
\end{lemma}

%%%%%%%%%%%%%%%%%%%%%%%%%%%%%%%%%%%%%%%%%%%%%%%%%%%%%%%%%%%%%%%%%%%%%%%%%%%%%%%%%%%%%%%%%%%%%%%%%%%%%%%%%%%%%%%%%%%%%%
%%%%%%%%%%%%%%%%%%%%%%%%%%%%%%%%%%%%%%%%%%%%%%%%%%%%%%%%%%%%%%%%%%%%%%%%%%%%%%%%%%%%%%%%%%%%%%%%%%%%%%%%%%%%%%%%%%%%%%
%%%%%%%%%%%%%%%%%%%%%%%%%%%%%%%%%%%%%%%%%%%%%%%%%%%%%%%%%%%%%%%%%%%%%%%%%%%%%%%%%%%%%%%%%%%%%%%%%%%%%%%%%%%%%%%%%%%%%%
Let us estimate the norm of $\left\|{\mathcal{S}^\epsilon_B}\right\|_{\mathcal{L}\left(L^2(\partial B), L^2(\partial B) \right)}$.
\begin{lemma}\label{lemmanrm1singulayer-casual}
  The operator norm of the compact operator $\mathcal{S}_{ D_\epsilon}:L^{2}(\partial D_\epsilon)\rightarrow L^{2}(\partial D_\epsilon)$, defined 
  in \eqref{defofSpartialDe}, is estimated by $\epsilon$,
i.e.\begin{eqnarray}\label{nrm1singulayer31-casual}
 \left\|\mathcal{S}_{ D_\epsilon}\right\|_{\mathcal{L}\left(L^2(\partial D_\epsilon), L^2(\partial D_\epsilon) \right)}
&\leq&\epsilon \left(\left\|{\mathcal{S}^{0}_B}\right\|_{\mathcal{L}\left(L^2(\partial B), L^2(\partial B) \right)}+\frac{1}{2\pi}\kappa\epsilon^2\vert\partial B\vert\right) ,
\end{eqnarray}
\end{lemma}
\begin{proofn}{\it{of Lemma \ref{lemmanrm1singulayer-casual}}.} To estimate the operator norm of $\mathcal{S}_{ D_\epsilon}$, we decompose 
$\mathcal{S}_{ D_\epsilon}=:\mathcal{S}_{ D_\epsilon}^\kappa=\mathcal{S}_{ D_\epsilon}^{i_{\kappa}}+\mathcal{S}_{ D_\epsilon}^{d_{\kappa}}$ into 
two parts $\mathcal{S}_{ D_\epsilon}^{i_{\kappa}}$ ( independent of $\kappa$ ) and  $\mathcal{S}_{ D_\epsilon}^{d_{\kappa}}$ ( dependent of $\kappa$ ) given by
\begin{eqnarray}
\mathcal{S}_{ D_\epsilon}^{i_{\kappa}} \psi(x):=\int_{\partial D_\epsilon} \frac{1}{4\pi|x-y|}\psi(y) dy,\label{def-of_Sik-normal}\\
\mathcal{S}_{ D_\epsilon}^{d_{\kappa}} \psi(x):=\int_{\partial D_\epsilon} \frac{e^{i\kappa|x-y|}-1}{4\pi|x-y|}\psi(y) dy.\label{def-of_Sdk-normal}
\end{eqnarray}
 With this definition, $\mathcal{S}_{ D_\epsilon}^{i_{\kappa}}:L^2(\partial D_\epsilon)\rightarrow L^2(\partial D_\epsilon)$ and $\mathcal{S}_{ D_\epsilon}^{d_{\kappa}}:L^2(\partial D_\epsilon)\rightarrow L^2(\partial D_\epsilon)$
 are compact. From \eqref{nrm1singulayer2dd}, it can be observed that
 \begin{equation}\label{nrm1singulayer2dd-bbb}
 \left\|\mathcal{S}_{ D_\epsilon}^{i_{\kappa}}\right\|_{\mathcal{L}\left(L^2(\partial D_\epsilon), L^2(\partial D_\epsilon) \right)}=
\epsilon\left\|{\mathcal{S}^{{i_\kappa}}_B}\right\|_{\mathcal{L}\left(L^2(\partial B), L^2(\partial B) \right)}=\epsilon\left\|{\mathcal{S}^{{0}}_B}\right\|_{\mathcal{L}\left(L^2(\partial B), L^2(\partial B) \right)}
\end{equation}
% with 
% \begin{eqnarray}
% \mathcal{S}^{\epsilon{i_\kappa}}_B \hat{\psi}(\xi):=\mathcal{S}^{{i_\kappa}}_B \hat{\psi}(\xi)=\int_{\partial B} \frac{1}{4\pi|\xi-\eta|}\hat{\psi}(\eta) d\eta,\label{def-of_Sik-casual}\\
% \mathcal{S}^{\epsilon{d_\kappa}}_B \hat{\psi}(\xi):=\int_{\partial B} \frac{e^{i\kappa\epsilon|\xi-\eta|}-1}{4\pi|\xi-\eta|}\hat{\psi}(\eta) d\eta.\label{def-of_Sdk-casual}
% \end{eqnarray}

On the other hand, as mentioned in (2.30) of  \cite[Lemma 2.5]{DPC-SM:MMS2013}, the following estimate can be obtained,
 \begin{eqnarray}\label{nrm1singulayer2dd-ccc}
 \left\|\mathcal{S}_{ D_\epsilon}^{d_{\kappa}}\right\|_{\mathcal{L}\left(L^2(\partial D_\epsilon), L^2(\partial D_\epsilon) \right)}
 &\leq&
\frac{1}{2\pi}\kappa\epsilon^3\vert\partial B\vert\, \mbox{ for }\kappa_{\max}diam(D_\epsilon)\leq 1.
\end{eqnarray}
Hence the result follows.
\end{proofn}

% \textcolor{red}{ Sir, here, I did not change the size from $(\frac{a}{2}+d^\alpha)$ to $(\frac{a}{2}+a^\alpha)$, as it seems the sufficient condition for the Neumann series looks 
% better when it involves $d$. Precisely, $\epsilon^{2-\beta}d^{-3\alpha}<c$ is better than $\epsilon^{2-\beta-3\alpha}<c$. So I feel we need to think in a different direction.} \\

\subsubsection{The multiple obstacle case}\label{gbmulobscase}
%\ ~ \ \\
\subsubsection*{A way of counting the small scatterers}\label{counting}
\begin{figure}[htp]
\centering
\includegraphics[width=4.5cm,height =4.5cm]{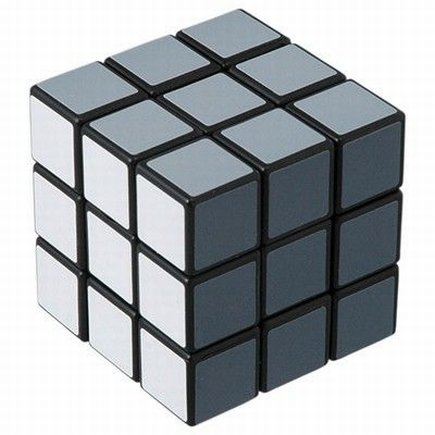}
\caption{Rubik's cube consisting of two layers}\label{fig:1-acsmall}
\end{figure}

Before proceeding further we make the following observation. For $m=1,\dots,M$ fixed, we distinguish between the obstacles $D_j$, $j\neq\,m$ by keeping them into different layers based on their distance from $D_m$.  Let $\Omega_m$, $1\leq\,m\leq\,M$ be the cubes
of center $z_m$ such that each side is of size $(\frac{a}{2}+d^\alpha)$ with $0\leq\alpha\leq{1}$ and it contains only $D_m$. Let us suppose that these cubes are arranged in a cuboid, for example unit rubics cube, see Fig \ref{fig:1-acsmall}, in different layers such that the total cubes upto the $n^{th}$
layer consists $(2n+1)^3$ cubes for $n=0,\dots,[d^{-\alpha}]$, and $\Omega_m$ is located on the center. Hence the number of obstacles 
located in the $n^{th}$, $n\neq0$ layer  will be $[(2n+1)^3-(2n-1)^3]$ and their distance from $D_m$ is more than ${n}d^\alpha$.
\bigskip

\bigskip

\begin{lemma} \label{propsmjsmmest}
 For $m,j=1,2,\dots,M$, the operator $\mathcal{S}_{mj}:L^2(\partial D_j)\rightarrow L^{2}(\partial D_m)$ defined 
%  by %in Proposition \ref{existence-of-sigmas} by
% \begin{eqnarray}\label{defofSmjed1}
%  \mathcal{S}_{mj}(\sigma_j)(s_m):=\int_{\partial D_j}\Phi_\kappa(s_m,s)\sigma_j(s)ds,\,s_m\in\partial D_m
% \end{eqnarray}
in Proposition \ref{existence-of-sigmas}, see \eqref{defofSmjed}, 
satisfies the following estimates,
\begin{itemize}
\item
For $j=m$,
\begin{eqnarray}\label{estinvsmm}
\left\|\mathcal{S}_{mm}\right\|_{\mathcal{L}\left(L^2(\partial D_m), L^2(\partial D_m) \right)}&\leq& \epsilon\left(\left\|\mathcal{S}^{{i_{\kappa}}}_{\c{B}}\right\|+\frac{1}{2\pi}\kappa\epsilon^2\vert\partial \c{B}\vert\right), %\label{{estinvsmm}}
\end{eqnarray}
for $\kappa_{\max} a\leq1$.
%where $C_{6m}:=\frac{2\pi\left\|{\mathcal{S}^{i_{\kappa}}_{B_m}}^{-1}\right\|_{\mathcal{L}\left(H^1(\partial B_m), L^2(\partial B_m) \right)}}{2\pi-(1+\kappa)\kappa\epsilon|\partial B_m|\left\|{\mathcal{S}^{i_{\kappa}}_{B_m}}^{-1}\right\|_{\mathcal{L}\left(H^1(\partial B_m), L^2(\partial B_m) \right)}}$.
\item For $j\neq m$, such that $D_j\in N^n_m, n=1,\dots,[d^{-\alpha}]$
\begin{eqnarray}\label{estinvsmj}
\left\|\mathcal{S}_{mj}\right\|_{\mathcal{L}\left(L^2(\partial D_j), L^2(\partial D_m)\right)}
                          &\leq&\frac{1}{4\pi}\frac{1}{nd^{\alpha}}\left|\partial \c{B} \right|\epsilon^{2},%\label{estinvsmj}
\end{eqnarray}
% \item For $j\neq m$, such that $D_j\in F_m$
% \begin{eqnarray}\label{estinvsmjF}
% \left\|\mathcal{S}_{mj}\right\|_{\mathcal{L}\left(L^2(\partial D_j), L^2(\partial D_m)\right)}
%                           &\leq&\frac{1}{4\pi}\frac{1}{d^\alpha}\left|\partial \c{B} \right|\epsilon^{2},%\label{estinvsmj}
% \end{eqnarray}
where $\left|\partial \c{B} \right|:=\max\limits_m \left|\partial B_m\right|$ and $\left\|\mathcal{S}^{{i_{\kappa}}}_{\c{B}}\right\|:=\max\limits_{m}\left\|\mathcal{S}^{{i_{\kappa}}}_{B_m}\right\|_{\mathcal{L}\left(L^2(\partial B_m), L^2(\partial B_m) \right)}$.

\end{itemize}
\end{lemma}
\begin{proofn}{\it{of Lemma \ref{propsmjsmmest}}.}
The estimate \eqref{estinvsmm} is nothing else but \eqref{nrm1singulayer31-casual} of Lemma \ref{lemmanrm1singulayer-casual}, replacing $B$ by $B_m$, $z$ by $z_m$ and $D_\epsilon$ by $D_m$ respectively. 
The proof of the estimate \eqref{estinvsmj} is a straightforward consequences of (2.37) in \cite[Lemma 2.6]{DPC-SM:MMS2013}. 
\end{proofn}
%%%%%%%%%%%%%%%%%%%%%%%%%%%%%%%%%%%%%%%%%%%%%%%%%%%%
\begin{proposition} \label{propsmjsmmestdbl-star}
 For $m,j=1,2,\dots,M$, the operator $\mathcal{D}_{mj}^*:L^2(\partial D_j)\rightarrow L^{2}(\partial D_m)$ defined 
in Proposition  \ref{existence-of-sigmas}, see \eqref{defofDmjed}, 
satisfies the following estimates,
\begin{itemize}
\item
For $j=m$,
\begin{eqnarray}\label{estinvsmmdbl-star}
 \left\|\left(-\frac{1}{2}I+\mathcal{D}_{mm}^*\right)^{-1}\right\|_{\mathcal{L}\left(L^2(\partial D_m), L^2(\partial D_m) \right)}&\leq&\grave{C}_{6m}, \label{invfnormD_{ii}1-star}
\end{eqnarray}
where $\grave{C}_{6m}:=\frac{2\pi\left\|\left(-\frac{1}{2}I+\mathcal{D}^{{i_{\kappa}}^*}_{B_m}\right)^{-1}\right\|_{\mathcal{L}\left(L^2(\partial B_m), L^2(\partial B_m) \right)}}
{2\pi-\kappa^2\epsilon^2|\partial B_m|\left\|\left(-\frac{1}{2}I+\mathcal{D}^{{i_{\kappa}}^*}_{B_m}\right)^{-1}\right\|_{\mathcal{L}\left(L^2(\partial B_m), L^2(\partial B_m) \right)}}$.
\item For $j\neq m$, such that $D_j\in N^n_m, n=1,\dots,[d^{-\alpha}]$
\begin{eqnarray}\label{estinvsmjdbl-star}
\left\|\mathcal{D}_{mj}^*\right\|_{\mathcal{L}\left(L^2(\partial D_j), L^2(\partial D_m)\right)}
                          &\leq&\frac{1}{4\pi}\left(\frac{\kappa}{nd^{\alpha}}+\frac{1}{n^2d^{2\alpha}}\right)\left|\partial \c{B} \right|\epsilon^{2},\label{fnormD_{ij}21-star}
\end{eqnarray}
where $\left|\partial \c{B} \right|:=\max\limits_m \partial B_m$.
\end{itemize}
In addition, as a consequence of \eqref{rep1sgldbllayer2}, we can also prove that 
\begin{itemize}
\item
For $j=m$,
\begin{eqnarray}\label{estinvsmmsgldbl-star}
 \left\|\left(-\frac{1}{2}I+\mathcal{D}_{mm}^*+\lambda_m \mathcal{S}_{mm}\right)^{-1}\right\|_{\mathcal{L}\left(L^2(\partial D_m), L^2(\partial D_m) \right)}&\leq&{C}_{6m}, \label{invfnormS_{ii}D_{ii}1-star}
\end{eqnarray}
where ${C}_{6m}:=\frac{2\pi\left\|\left(-\frac{1}{2}I+\mathcal{D}^{{i_{\kappa}}^*}_{B_m}+\lambda_m\epsilon \mathcal{S}^{{i_{\kappa}}}_{B_m}\right)^{-1}\right\|_{\mathcal{L}\left(L^2(\partial B_m), L^2(\partial B_m) \right)}}
{2\pi-(\vert\lambda_m\vert+\kappa)\kappa\epsilon^2|\partial B_m|\left\|\left(-\frac{1}{2}I+\mathcal{D}^{{i_{\kappa}}^*}_{B_m}+\lambda_m\epsilon \mathcal{S}^{{i_{\kappa}}}_{B_m}\right)^{-1}\right\|_{\mathcal{L}\left(L^2(\partial B_m), L^2(\partial B_m) \right)}}$.
\end{itemize}
\end{proposition}
\begin{proofn}{\it{of Proposition \ref{propsmjsmmestdbl-star}}.}
This result can be proved in the similar lines of the proof of \cite[Proposition 2.17]{DPC-SM:MMS2013}.
\end{proofn}

%%%%%%%%%%%%%%%%%%%%%%%%%%%%%%%%%%%%%%%%%%%%%%%%%%%%
\begin{proofe}
\textbf{\textit{End of the proof of Proposition \ref{normofsigmastmt}}.}
 By substituting  (\ref{estinvsmj}) and (\ref{fnormD_{ij}21-star}) in \eqref{DKSnrm}, 
 \eqref{invfnormS_{ii}D_{ii}1-star} in \eqref{invDLSnrm}, and using our discusion 
 related to Fig \ref{fig:1-acsmall} on how we count the number of small scatterers, 
  we obtain
\begin{eqnarray}
\left\|DK^*+\lambda K\right\|
    &\equiv&%\sum\limits_{m=1}^{M}\sum\limits_{\substack{j=1\\j\neq m}}^{M}
\max\limits_{1\leq m \leq M}\sum_{\substack{j=1\\j\neq\,m}}^{M}\left\|\mathcal{D}_{mj}^*+\lambda_m \mathcal{S}_{mj}\right\|_{\mathcal{L}\left(L^{2}(\partial D_j),L^{2}(\partial D_m)\right)}\nonumber\\
   &\leq&\sum_{n=1}^{[d^{-\alpha}]}[(2n+1)^3-(2n-1)^3]\frac{1}{4\pi}\left(\frac{\vert\lambda\vert+\kappa}{nd^{\alpha}}+\frac{1}{n^2d^{2\alpha}}\right)\left|\partial \c{B} \right|\epsilon^{2}\nonumber\\
   &=&\sum_{n=1}^{[d^{-\alpha}]}[24n^2+2]\frac{1}{4\pi}\left(\frac{\vert\lambda\vert+\kappa}{nd^{\alpha}}+\frac{1}{n^2d^{2\alpha}}\right)\left|\partial \c{B} \right|\epsilon^{2}\nonumber\\
   &=&\sum_{n=1}^{[d^{-\alpha}]}\frac{1}{2\pi}[12n^2+1]\left(\frac{\vert\lambda\vert+\kappa}{nd^{\alpha}}+\frac{1}{n^2d^{2\alpha}}\right)\left|\partial \c{B} \right|\epsilon^{2}\nonumber\\
   &=&\frac{1}{2\pi}\left[(\vert\lambda\vert+\kappa)d^{-\alpha}\sum_{n=1}^{[d^{-\alpha}]}[12n+\frac{1}{n}]+d^{-2\alpha}\sum_{n=1}^{[d^{-\alpha}]}[12+\frac{1}{n^2}]\right]\left|\partial \c{B} \right|\epsilon^{2}\nonumber\\
    &\leq&\frac{1}{2\pi}\left[(\vert\lambda\vert+\kappa)(6d^{-3\alpha}+7d^{-2\alpha})+13d^{-3\alpha}\right]\left|\partial \c{B} \right|\epsilon^{2}\label{KDKnrm1}
\end{eqnarray}
and
\begin{eqnarray}
\left\|\left(-\frac{1}{2}\textbf{I}+DL^*+\lambda L \right)^{-1}\right\|
     &\equiv&\max\limits_{1\leq m \leq M}\left\|\left(-\frac{1}{2}I+\mathcal{D}_{mm}^*+\lambda_m \mathcal{S}_{mm}\right)^{-1}\right\|_{\mathcal{L}\left(L^{2}(\partial D_m),L^{2}(\partial D_m)\right)}\nonumber\\ 
     &\equiv&\max\limits_{1\leq m \leq M}{C}_{6m}.\label{invLDLnrm1}
\end{eqnarray}

Hence, (\ref{KDKnrm1}-\ref{invLDLnrm1}) provides
\begin{eqnarray}
\left\|(-\frac{1}{2}\textbf{I}+DL^{*}+\lambda L)^{-1}\right\|\left\|DK^{*}+\lambda K\right\|\nonumber\\
    &\hspace{-2.5cm}\leq&\hspace{-1cm}\underbrace{\left(\max\limits_{1\leq m \leq M}{C}_{6m}\right)\left|\partial \c{B}\right|\frac{1}{2\pi}\left[(\vert\lambda\vert+\kappa)(6+7d^{\alpha})+13\right]d^{-3\alpha}\epsilon^{2}}_{=:{C}_s}.\label{invDLDKKnrm1}
\end{eqnarray}
By imposing the condition $\left\|(-\frac{1}{2}\textbf{I}+DL^{*}+\lambda L)^{-1}\right\|\left\|DK^{*}+\lambda K\right\|<1$, we have from \eqref{nrminvLplusK} and (\ref{sigmaU^Inrm1}-\ref{sigmaU^Inrm});
\begin{eqnarray}\label{nrminvDLplusDK2}
 \left\|\sigma_m\right\|_{L^{2}(\partial D_m)}\leq\left\|\sigma\right\| 
                                  &\leq&\frac{\left\|(-\frac{1}{2}\textbf{I}+DL^{*}+\lambda L)^{-1}\right\|}{1-\left\|(-\frac{1}{2}\textbf{I}+DL^{*}+\lambda L)^{-1}\right\|\left\|DK^{*}+\lambda K\right\|}\left\|(\partial_\nu+\lambda)U^{In}\right\|\nonumber\\
                   &\leq&{C}_p\left\|(-\frac{1}{2}\textbf{I}+DL^{*}+\lambda L)^{-1}\right\| \max\limits_{1\leq m \leq M}\left\|(\partial_\nu+\lambda)U^{i}\right\|_{L^{2}(\partial D_m)}\hspace{.25cm} \left( {C}_p\geq\frac{1}{1-{C}_s}\right)\nonumber\\
                   &\substack{\leq\\ \eqref{invLDLnrm1} }&\mathrm{C}  \max\limits_{1\leq m \leq M}\left\|(\partial_\nu+\lambda)U^{i}\right\|_{L^{2}(\partial D_m)}\hspace{.25cm} \left(\mathrm{C}:={C}_p\max\limits_{1\leq m \leq M}{C}_{6m}\right),
\end{eqnarray}
for all $m\in\{1,2,\dots,M\}$.
But,
\begin{eqnarray}\label{normU^i}
\left\|(\partial_\nu+\lambda)U^{i}\right\|_{L^{2}(\partial D_m)}&\leq&\vert\lambda\vert\left\|U^{i}\right\|_{L^{2}(\partial D_m)}+\left\|\partial_\nu U^{i}\right\|_{L^{2}(\partial D_m)}\nonumber\\
                                            &=&\vert\lambda\vert\epsilon~\left|\partial B_m\right|^{\frac{1}{2}}+k\epsilon\left|\partial B_m\right|^{\frac{1}{2}}\quad\left(\mbox{Since },U^{i}(x,\theta)=e^{i\kappa{x}\cdot\theta}\right) \nonumber\\
                                            &\leq&(\vert\lambda\vert+\kappa)\epsilon\left|\partial \c{B}\right|^{\frac{1}{2}}, \forall m=1,2,\dots,M.
\end{eqnarray}
Now by substituting \eqref{normU^i} in \eqref{nrminvDLplusDK2}, for each $m=1,\dots,M$, we obtain
\begin{eqnarray}\label{nrmsigmaf}
\left\|\sigma_m\right\|_{L^{2}(\partial D_m)}&\,\leq\,&  {\mathcal{C}}(\kappa)\epsilon,
\end{eqnarray}
where $\hspace{.25cm}\mathcal{C}(\kappa):=\mathrm{C} \left|\partial \c{B}\right|^{\frac{1}{2}}(\vert\lambda\vert+\kappa)$. 
\par The condition $\left\|(-\frac{1}{2}\textbf{I}+DL^{*}+\lambda L)^{-1}\right\|\left\|DK^{*}+\lambda K)\right\|<1$ is satisfied if 
\begin{eqnarray}\label{Accond-invLK-singl-small}
{C}_s&=&\left(\max\limits_{1\leq m \leq M}{C}_{6m}\right)\left|\partial \c{B}\right|\frac{1}{2\pi}\left[(\vert\lambda\vert+\kappa)(6+7d^{\alpha})+13\right]d^{-3\alpha}\epsilon^{2}\,<\,1.
\end{eqnarray}
Since $\lambda_m=\lambda_{m0}\, a^{-\beta}$ and $d\approx a^t$, in particular $d \geq d_{min} a^t$, then 
%\textcolor{red}{%Since $(2\kappa+1)d<\bar{c}$ for $\bar{c}:=(2\kappa_{\max}+1)d_{\max}$, then \eqref{Accond-invLK-singl-small} 
\eqref{Accond-invLK-singl-small}
reads as $a^{-3\alpha t+2-\beta}<\grave{c}$, where we set 
\begin{equation}\label{lambda-}
\grave{c}:=\left(\big[(\lambda_+ +\kappa_{\max}a^{\beta})(6+7{[d_{\max}]}^{\alpha})+
13a^{\beta}\big]\frac{1}{2\pi}\frac{|\partial\c{B}|}{[\max\limits_{1\leq m \leq M}diam(B_m)]^2}\max\limits_{1\leq m \leq M}C_{6m}d^{-\alpha}_{min}\right)^{-1}>1\,
\end{equation}
with $\lambda_+:=\max\limits_{1\leq m \leq M}\vert\lambda_{m0}\vert$ and with the rewritten form of ${C}_{6m}$ mentioned in Proposition \ref{propsmjsmmestdbl-star} as 
$$
{C}_{6m}:=\frac{2\pi\left\|\left(-\frac{1}{2}I+\mathcal{D}^{{i_{\kappa}}^*}_{B_m}+\frac{\lambda_{m0}\epsilon^{1-\beta}}{[\max\limits_{1\leq m \leq M} diam(B_m)]^\beta} 
\mathcal{S}^{{i_{\kappa}}}_{B_m}\right)^{-1}\right\|_{\mathcal{L}\left(L^2(\partial B_m), L^2(\partial B_m) \right)}}
{2\pi-(\frac{\lambda_{m0}\kappa\epsilon^{2-\beta}}{[\max\limits_{1\leq m \leq M}diam(B_m)]^\beta}+\kappa^2\epsilon^{2})
|\partial B_m|\left\|\left(-\frac{1}{2}I+\mathcal{D}^{{i_{\kappa}}^*}_{B_m}+\frac{\lambda_{m0}\epsilon^{1-\beta}}{[\max\limits_{1\leq m \leq M}diam(B_m)]^\beta} 
\mathcal{S}^{{i_{\kappa}}}_{B_m}\right)^{-1}\right\|_{\mathcal{L}\left(L^2(\partial B_m), L^2(\partial B_m) \right)}}.
$$
Observe that $s=3 \alpha t$. Hence \eqref{Accond-invLK-singl-small} makes sense if $s\leq 2-\beta$ and $\lambda_+$ satisfies (\ref{lambda-}).

Again, since $\lambda_m=\lambda_{m0} a^{-\beta}$, \eqref{nrmsigmaf} can be rewritten as 
\begin{eqnarray}\label{nrmsigmaf-1}
\left\|\sigma_m\right\|_{L^{2}(\partial D_m)}&\,\leq\,&  \grave{\mathcal{C}}(\kappa)\epsilon^{1-\beta},
\end{eqnarray}
with\footnote{It is important to remark that, if we do not distinguish the 
near by and far obstacles, as it is discussed in the beginning of section
\ref{gbmulobscase}, and by following the way it was done in \cite{DPC-SM:MMS2013,DPC-SM:MANA2015} we can get the estimate
$\left\|DK^*+\lambda K\right\|\leq\frac{M-1}{4\pi}\left(\frac{\vert\lambda\vert+\kappa}{d}+\frac{1}{d^2}\right)\left|\partial \c{B} \right|\epsilon^{2}$ in place of \eqref{KDKnrm1} and 
hence the %second condition of \eqref{conditions} in Theorem \ref{Maintheorem-ac-small-sing} will be replaced by
condition \eqref{Accond-invLK-singl-small} will be replaced by
$(M-1)\frac{a^{2-\beta}}{d^2}<c_0$, for some suitable constant $c_0$. 
%Infact it is a sufficient condition for $C_s<1$ with $\grave{\grave{c}}:=\left(\frac{(\vert\lambda_0\vert+\kappa_{\max})d_{\max}+1}{4\pi}|\partial\c{B}|\max\limits_{1\leq m \leq M}C_{6m}\right)^{-1}$. 
However, this condition is too strong to enable us to apply our asymptotic expansion to the effective medium theory where we need to choose $M \sim a^{-s}$ with
$ s=2-\beta$ and $d \sim a^t$with $t\geq \frac{s}{3}$.
} 
$$\hspace{.25cm}\grave{\mathcal{C}}(\kappa):=\mathrm{C} \left|\partial \c{B}\right|^{\frac{1}{2}}(\frac{\lambda_+}{[\max\limits_{1\leq m \leq M}diam(B_m)]^\beta}+\kappa\epsilon^{\beta})\,\left(<\mathrm{C} \left|\partial \c{B}\right|^{\frac{1}{2}}(\frac{\lambda_+}{[\max\limits_{1\leq m \leq M} diam(B_m)]^\beta}+\kappa_{\max})\right).$$

\end{proofe}

%\subsection{Approximation of the far-fields}\label{SLPR-3}
\subsection{Approximation of the far-fields. I. Approximation by 
the total charges}\label{SLPR-3}

We start with the definition of the total charges $Q_m,\,m=1,\dots\,M$.
\bigskip

\begin{definition}%[Definition \ref{Qmdef}.]
\label{Qmdef}
We call the $\sigma_m$'s used in \eqref{qcimprequiredfrm1}, the solution of the problem (\ref{acimpoenetrable}-\ref{radiationc}), the surface charge distributions. 
Using these surface charge distributions, we define the total charge on each surface $\partial D_m$ denoted by $Q_m$ as
\begin{eqnarray}\label{defofQm}
 Q_m:=\int_{ \partial D_m} \sigma_m(s) ds.
\end{eqnarray}
\end{definition}

In the following proposition, we provide
an approximate of the far-fields in terms of the total charges $Q_m$.
\bigskip

\begin{proposition}\label{farfldthm} 
The far-field pattern $U^\infty$ of the scattered solution of the problem 
(\ref{acimpoenetrable}-\ref{radiationc}) has the following asymptotic
expansion
\begin{equation}\label{x oustdie1 D_m}
%U^t(x)=U^i(x)+\sum_{m=1}^{M}\Phi_\kappa(x,z_{m})\left[Q_m+O\left(\kappa a+\frac{a}{\textcolor{red}{|x-z_m|}}\right)\int_{\partial D_m} \vert\sigma_{m}(s)\vert ds\right],
U^\infty(\hat{x})=\sum_{m=1}^{M}[e^{-i\kappa\hat{x}\cdot z_{m}}Q_m+O(\kappa\,a^{3-\beta})],
\end{equation}
with $Q_m$ given by (\ref{defofQm}), if 
$\kappa_{\max}\,a<1$ where $O(\kappa\,a^{3-\beta})\,\leq\,
C\kappa\,a^{3-\beta}$ and $C:=\frac{|\partial\c{B}|\mathrm{C} 
(\lambda_+ +\kappa_{\max})}{\left(\max\limits_{1\leq m \leq M} 
diam(B_m)\right)^{2-\beta}}$.%\hspace{.25cm} \textcolor{violet}{(\mbox{Since}, a=\epsilon\max\limits_{1\leq m \leq M} diam(B_m).)}

\end{proposition}
\begin{proofn}{\it{of Proposition \ref{farfldthm}}.}
From \eqref{qcimprequiredfrm1}, we have 
\begin{eqnarray*}
 U^{s}(x)&=&\sum_{m=1}^{M}\int_{\partial D_m}\Phi_\kappa(x,s)\sigma_{m} (s)ds,\text{ for }x\in\mathbb{R}^{3}\backslash\left(\mathop{\cup}\limits_{m=1}^M \bar{D}_m\right). \nonumber
\end{eqnarray*}
Hence
\begin{eqnarray}\label{xfarawayimpnt}
U^{\infty}(\hat{x})&=&\sum_{m=1}^{M}\int_{\partial D_m}e^{-i\kappa\hat{x}\cdot s}\sigma_{m} (s)ds\nonumber\\
 &=&\sum_{m=1}^{M}\left(\int_{\partial D_m}e^{-i\kappa\hat{x}\cdot\,z_{m}}\sigma_{m}(s)ds+\int_{\partial D_m}[e^{-i\kappa\hat{x}\cdot\,s}-e^{-i\kappa\hat{x}\cdot\,z_{m}}]\sigma_{m}(s)ds\right)\nonumber\\
 &=&\sum_{m=1}^{M}\left(e^{-i\kappa\hat{x}\cdot\,z_{m}}Q_m+\int_{\partial D_m}[e^{-i\kappa\hat{x}\cdot\,s}-e^{-i\kappa\hat{x}\cdot\,z_{m}}]\sigma_{m}(s)ds\right).
%&=&\sum_{m=1}^{M}e^{-i\kappa\hat{x}\cdot\,z_{m}}Q_m+O(\kappa\,a^2).
\end{eqnarray}
  %For every $m=1,2,\dots,M$,  we have from Proposition \ref{normofsigmastmt};
For every $m=1,2,\dots,M$, we have from Proposition \ref{normofsigmastmt};
% \begin{eqnarray}\label{estimationofintsigma}
% \int_{\partial D_m}\vert\sigma_{m}(s)\vert ds
% &\leq&Ca^{2-\beta},
% \end{eqnarray}
\begin{eqnarray}\label{estimationofintsigma}
 \left|\int_{\partial D_m}\vert\sigma_{m}(s)\vert ds\right|%&=&\left|\int_{ \partial D_m} \sigma_m(s) ds\right|\nonumber\\
    &\leq& \| 1  \|_{L^2(\partial D_m)} \| \sigma_m  \|_{L^2(\partial D_m)}\nonumber\\
    &\substack{\leq \\ \eqref{nrmsigmaf-1}}& \| 1  \|_{L^2(\partial D_m)}\mathrm{C} \left|\partial \c{B}\right|^{\frac{1}{2}}(\lambda_+ +\kappa_{\max})\,\epsilon^{1-\beta} \nonumber\\
    &\leq&|\partial \c{B}|\mathrm{C} (\frac{\lambda_+}{[\max\limits_{1\leq m \leq M} diam(B_m)]^{\beta}}+\kappa_{\max})\,\epsilon^{2-\beta}
\end{eqnarray}
\begin{equation}
 \text{with}\quad C:=%\frac{\tilde{c}}{\left(\max\limits_{1\leq m \leq M} diam(B_m)\right)^{2-\beta}}
 =\frac{|\partial \c{B}|\mathrm{C} (\lambda_+ +\kappa_{\max}[\max\limits_{1\leq m \leq M} diam(B_m)]^{\beta})}{\left(\max\limits_{1\leq m \leq M} diam(B_m)\right)^{2}}.%\hspace{.25cm} \textcolor{violet}{(\mbox{Since}, a=\epsilon\max\limits_{1\leq m \leq M} diam(B_m).)}
\end{equation}
It gives us the following estimate;
\begin{eqnarray}\label{acestimateforexponentsdif-pre}
\left|\int_{\partial D_m}[e^{-i\kappa\hat{x}\cdot\,s}-e^{-i\kappa\hat{x}\cdot\,z_{m}}]\sigma_{m}(s)ds\right|
                       %  &\leq&\int_{\partial D_m}\left|e^{-i\kappa\hat{x}\cdot\,s}-e^{-i\kappa\hat{x}\cdot\,z_{m}}\right \|\sigma_{m}(s)|ds\nonumber\\
                         &\leq&\int_{\partial D_m}\left|e^{-i\kappa\hat{x}\cdot\,s}-e^{-i\kappa\hat{x}\cdot\,z_{m}}\right| |\sigma_{m}(s)|ds\nonumber\\
                         &\leq&\int_{\partial D_m}\sum_{l=1}^{\infty}\kappa^l|s-z_{m}|^l|\sigma_{m}(s)|ds\nonumber\\
                         &\leq&\int_{\partial D_m}\sum_{l=1}^{\infty}\kappa^l\left(\frac{a}{2}\right)^l|\sigma_{m}(s)|ds\nonumber\\
                        % &\leq&\sum_{l=1}^{\infty}\kappa^la^l\int_{\partial D_m}|\sigma_{m}(s)|ds\nonumber\\
                         &\substack{\leq\\ \eqref{estimationofintsigma}}&Ca^{2-\beta}\sum_{l=1}^{\infty}\kappa^l\left(\frac{a}{2}\right)^l\nonumber\\
                         &=&\frac{1}{2}C\kappa\,a^{3-\beta}\frac{1}{1-\kappa\frac{a}{2}},\text{ if }
a<\frac{2}{\kappa_{\max}}\left(\leq\frac{2}{\kappa}\right)
\end{eqnarray}
which means
\begin{eqnarray}\label{acestimateforexponentsdif}
\int_{\partial D_m}[e^{-i\kappa\hat{x}\cdot\,s}-e^{-i\kappa\hat{x}\cdot\,z_{m}}]\sigma_{m}(s)ds&\leq&C\kappa\,a^{3-\beta},\,\text{for}\,a\leq\frac{1}{\kappa_{\max}}.
\end{eqnarray}
Now substitution of \eqref{acestimateforexponentsdif} in \eqref{xfarawayimpnt} gives the required result \eqref{x oustdie1 D_m}.
\end{proofn}

\subsection{Approximation of the far-fields. II. Estimates of the total charges}
\subsubsection{Derivation of the linear algebraic system}
We start with following a priori estimate on the total charges $Q_m, \; m=1,...,M$.
\begin{lemma}\label{Qmestbigo}
For $m=1,2,\dots,M$, we have%the absolute value of the total charge $Q_m$ on each surface $\partial D_m$ is bounded by $\epsilon$. i.e. 
\begin{eqnarray}\label{estofQm}
 |Q_m|&\leq&\tilde{c}\,\epsilon^{2-\beta},
\end{eqnarray}
where $\tilde{c}:=|\partial \c{B}|\mathrm{C} (\frac{\lambda_+}{[\max\limits_{1\leq m \leq M} diam(B_m)]^{\beta}}+\kappa_{\max})$ with $\partial\c{B}$ and $\mathrm{C}$ are defined in \eqref{estinvsmj} and \eqref{nrminvDLplusDK2} respectively.
\end{lemma}
\begin{proofn}{\it{of Lemma \ref{Qmestbigo}}.}
% From Proposition \ref{normofsigmastmt}, we know that the surface charge distributions $\sigma_m\in L^2(\partial D_m)$ are bounded by a constant, $\grave{\mathcal{C}}(\kappa)\,\left(<\mathrm{C} \left|\partial \c{B}\right|^{\frac{1}{2}}(\lambda_0+\kappa_{\max})\right)$, times $\epsilon^{1-\beta}$ . Hence
The proof follows as below;
\begin{eqnarray*}\label{Q2}
 |Q_m|&=&\left|\int_{ \partial D_m} \sigma_m(s) ds\right|\nonumber\\
    &\leq& \| 1  \|_{L^2(\partial D_m)} \| \sigma_m  \|_{L^2(\partial D_m)}\nonumber\\
    &\substack{\leq\\ \eqref{nrmsigmaf-1}}& \| 1  \|_{L^2(\partial D_m)}\mathrm{C} \left|\partial \c{B}\right|^{\frac{1}{2}}(\frac{\lambda_+}{[\max\limits_{1\leq m \leq M} diam(B_m)]^{\beta}}+\kappa_{\max})\,\epsilon^{1-\beta} \nonumber\\
    &\leq&|\partial \c{B}|\mathrm{C} (\frac{\lambda_+}{[\max\limits_{1\leq m \leq M} diam(B_m)]^{\beta}}+\kappa_{\max})\,\epsilon^{2-\beta}.
\end{eqnarray*}
% Hence, we obtain that
% $$Q_m=O(\epsilon+\textcolor{red}{\kappa}\epsilon).$$
\end{proofn}
The following proposition gives an approximate characterization of the total charges 
$Q_m,\; m=1,...,M$.
\begin{proposition} \label{fracqfracc-ac}For $m=1,2,\dots,M$, the total charge ${Q}_m$ on each surface $\partial D_m$ of the small scatterer $D_m$ can be calculated from the algebraic system 
  \begin{eqnarray}\label{Q_mintdbl-4}
 \frac{{Q}_m}{\bar{C}_m} &=&-U^{i}(z_m)-\sum_{\substack{j=1 \\ j\neq m}}^{M}\bar{C}_j \Phi_\kappa(z_m,z_j)\frac{{Q}_j}{\bar{C}_j}+ 
Err
  \end{eqnarray}
  where $Err:=O\left( a^{1-\beta}+\frac{a^{3-\beta}}{d^{3\alpha}}\right)$  and $\bar{C}_m := -\lambda_m|\partial D_m|$.

\end{proposition}
\begin{proofn}{\it{of Proposition \ref{fracqfracc-ac}}.}
% Let us derive a formula for $Q_m$. 
For $s_m\in \partial D_m$, using the impedance boundary condition \eqref{acgoverningsupport}, we have
\begin{eqnarray}\label{Q_mintdbl}
 0&=&\frac{\partial U^{t}}{\partial \nu_m}(s_m)+\lambda_m U^{t} (s_m)=
-\frac{\sigma_{m} (s_m)}{2}+\int_{\partial D_m}\frac{\partial\Phi_\kappa}{\partial \nu_m}(s_m,s)\sigma_{m} (s)ds+\sum_{\substack{j=1\\j\neq m}}^{M}\int_{\partial D_j}\frac{\partial\Phi_\kappa}{\partial \nu_m}(s_m,s)\sigma_{j} (s)ds\nonumber\\
&& +\lambda_m\sum_{j=1}^{M}\int_{\partial D_j}\Phi_\kappa(s_m,s)\sigma_{m} (s)ds+\frac{\partial U^i}{\partial \nu_m}{(s_m)}+\lambda_m U^i{(s_m)}\nonumber\\
\end{eqnarray}
Integrating the above on $\partial D_m$, we can write it as
\begin{equation}%\label{Q_mintdbl-1}
\begin{split}
-\frac{1}{2}\int_{\partial D_m}\sigma_{m} (s_m) d{s_m}+\int_{\partial D_m}\left(\int_{\partial D_m}\frac{\partial\Phi_\kappa}{\partial \nu_m}(s_m,s) d{s_m}\right) \sigma_{m} (s)ds &+\sum_{\substack{j=1\\j\neq m}}^{M}\int_{\partial D_j}\left(\int_{\partial D_m}\frac{\partial\Phi_\kappa}{\partial \nu_m}(s_m,s)d{s_m}\right) \sigma_{j} (s)ds
 \nonumber\\
 +\lambda_m\sum_{j=1}^{M}\int_{\partial D_j}\left(\int_{\partial D_m}\Phi_\kappa(s_m,s) d{s_m}\right) \sigma_{m} (s)ds&=-\int_{\partial D_m}\frac{\partial U^i}{\partial \nu_m}{(s_m)} d{s_m}-\int_{\partial D_m}\lambda_m U^i{(s_m)}d{s_m}\nonumber\\
\end{split}
\end{equation}
It can be rewritten as
\begin{equation}\label{Q_mintdbl-1}
\begin{split}
-\frac{1}{2}Q_m&+\underbrace{\int_{\partial D_m}\left(\int_{\partial D_m}\frac{\partial\Phi_{0}}{\partial \nu_m}(s_m,s) d{s_m}\right) \sigma_{m} (s)ds}_{=: A} +
\underbrace{\sum_{\substack{j\neq m}}^{M}\int_{\partial D_j}\left(\int_{\partial D_m}\frac{\partial\Phi_\kappa}{\partial \nu_m}(s_m,s)d{s_m}\right) \sigma_{j} (s)ds}_{=: B}
 \nonumber\\
 &+\lambda_m\underbrace{\int_{\partial D_m}\left(\int_{\partial D_m}\Phi_\kappa(s_m,s) d{s_m}\right) \sigma_{m} (s)ds}_{=: C}+\lambda_m\underbrace{\sum_{j\neq m}^{M}\int_{\partial D_j}\left(\int_{\partial D_m}\Phi_\kappa(s_m,z_j) d{s_m}\right) \sigma_{j} (s)ds}_{=: D}
 \nonumber\\
 &=-\int_{\partial D_m}\frac{\partial U^i}{\partial \nu_m}{(s_m)} d{s_m}-\int_{\partial D_m}\lambda_m U^i{(s_m)}d{s_m}+A^{\prime}+\lambda_m D^{\prime},%\nonumber\\
\end{split}
\end{equation}
with
\begin{eqnarray}
A{^\prime}&:=&\int_{\partial D_m}\left(\int_{\partial D_m}\left[\frac{\partial\Phi_\kappa}{\partial \nu_m}(s_m,s)-\frac{\partial\Phi_0}{\partial \nu_m}(s_m,s)\right] d{s_m}\right) \sigma_{m} (s)ds \label{defofAprm}\\
D{^\prime}&:=&\sum_{j\neq m}^{M}\int_{\partial D_j}\left(\int_{\partial D_m}\left[\Phi_\kappa(s_m,s)-\Phi_\kappa(s_m,z_j)\right] d{s_m}\right) \sigma_{j} (s)ds \label{defofDprm}. 
\end{eqnarray}
% \clearpage
% % % % For $z\in\bar{D}_j, j\neq m$, let us consider the integral terms $\int_{\partial D_m}\frac{\partial\Phi_\kappa}{\partial \nu_m}(s_m,z)d{s_m}$ and $\int_{\partial D_m}\Phi_\kappa(s_m,z)d{s_m}$;
\begin{itemize}
\item

We can approximate $A$ and $A^{\prime}$ as follows;
\begin{eqnarray}\label{approximating-A}
A&=&\int_{\partial D_m}\left(\int_{\partial D_m}\frac{\partial\Phi_0}{\partial \nu_m}(s_m,s) d{s_m}\right) \sigma_{m} (s)ds\nonumber\\
&=&\int_{\partial D_m}[\mathbf{K}^{i_{\kappa}}_{ D_m}(1)](s) \sigma_{m} (s)ds\nonumber\\
&=&-\frac{1}{2}\int_{\partial D_m} \sigma_{m} (s)ds\nonumber\\
&=&-\frac{1}{2}Q_m.
\end{eqnarray}
Here $\mathbf{K}^{i_{\kappa}}_{ D_m}$ is the double layer operator defined as in \eqref{nbc23} but with zero frequency and on the boundary of $D_m$. Observe that,

\begin{eqnarray}\label{approximating-Apr}
|A{^\prime}|&=&\left\vert\int_{\partial D_m}\left(\int_{\partial D_m}\left[\frac{\partial\Phi_\kappa}{\partial \nu_m}(s_m,s)-\frac{\partial\Phi_0}{\partial \nu_m}(s_m,s)\right] d{s_m}\right) \sigma_{m} (s)ds\right\vert\nonumber\\
&=&\left\vert\int_{\partial D_m}\left(\int_{\partial D_m}\left[\frac{i\kappa(s_m-s)\cdot \nu_m(s_m)}{4\pi|s_m-s|^2}\sum_{l=1}^\infty (i\kappa|s_m-s|)^l\left(\frac{1}{l!}-\frac{1}{(l+1)!}\right)\right] d{s_m}\right) \sigma_{m} (s)ds\right\vert\nonumber\\
&\leq&\left\vert\int_{\partial D_m}\left(\int_{\partial D_m}\left[\frac{\kappa^2}{4\pi}\sum_{l=0}^\infty \frac{(\kappa{a})^l}{2^l}\right] d{s_m}\right) |\sigma_{m} (s)|ds\right\vert\nonumber\\
&=&O\left(\kappa^2 a^{4-\beta}\right).%\textcolor{red}{\rightarrow \mbox{ this term may not be required, if I can show that   $\mathcal{K}^*{1}=\mathcal{K}{1}=-1$ on $\partial D$ even for non zero frequencies.}}
\end{eqnarray}
\item Now, we can approximate $B$
\begin{eqnarray}\label{approximating-B}
B&=&\sum_{\substack{j\neq m}}^{M}\int_{\partial D_j}\left(\int_{\partial D_m}\frac{\partial\Phi_\kappa}{\partial \nu_m}(s_m,s)d{s_m}\right) \sigma_{j} (s)ds\nonumber\\
&=&\sum_{\substack{j\neq m}}^{M}\int_{\partial D_j}\left(\int_{ D_m}\Delta\Phi_\kappa(y_m,s) d{y_m}\right) \sigma_{j} (s)ds\nonumber\\
&=&\sum_{\substack{j\neq m}}^{M}\int_{\partial D_j}\left(\int_{ D_m}[-\kappa^2\Phi_\kappa(y_m,s)] d{y_m}\right) \sigma_{j} (s)ds\nonumber\\
 &=& \sum_{\substack{j\neq m}}^{M}\int_{\partial D_j}\left[O\left(\kappa^2\frac{a^3}{d_{mj}}\right)\right]\sigma_{j} (s)ds\nonumber\\
 &=& \sum_{n=1}^{[d^{-\alpha}]}[(2n+1)^3-(2n-1)^3]\left[O\left(\kappa^2\frac{a^{5-\beta}}{nd^{\alpha}}\right)\right]\nonumber\\
  &=& O\left(2\kappa^2\frac{a^{5-\beta}}{d^{\alpha}}\sum_{n=1}^{[d^{-\alpha}]}[12n+\frac{1}{n}]\right)\nonumber\\
  &=& O\left(2\kappa^2\frac{a^{5-\beta}}{d^{2\alpha}}[\frac{6}{d^{\alpha}}+7]\right).
\end{eqnarray}

\item Since $\Phi_\kappa(s_m,s)-\Phi_0(s_m,z_j)=O\left(\frac{\kappa{a}}{d_{mj}}+\frac{a}{d_{mj}^2}\right)$ for $s\in \partial{D}_j, j\neq m$, we can approximate $D$ and $D^{\prime}$ as follows;
\begin{eqnarray}\label{approximating-D}
D&=&\sum_{\substack{j\neq m}}^{M}\int_{\partial D_j}\left(\int_{\partial D_m}\Phi_\kappa(s_m,z_j) d{s_m}\right) \sigma_{j} (s)ds
\nonumber\\&=&\sum_{\substack{j\neq m}}^{M}\left[
\int_{\partial D_j}\left(\int_{\partial D_m}\Phi_\kappa(z_m,z_j) d{s_m}\right) \sigma_{j} (s)ds\right.\nonumber\\
&&\left.\qquad+\int_{\partial D_j}\left(\int_{\partial D_m}[\Phi_\kappa(s_m,z_j)-\Phi_\kappa(z_m,z_j)] d{s_m}\right) \sigma_{j} (s)ds\right]\nonumber\\
% % %  &=&\textcolor{blue}{\sum_{\substack{j\neq m}}^{M}\left[\Phi_\kappa(z_m,z_j) Q_j |\partial D_m|+O\left(\frac{a}{d}\left(\kappa+\frac{1}{d}\right)|\partial D_m| \left|\int_{\partial D_j} \sigma_{j} (s)ds\right|\right)\right]}\nonumber\\
% % %  &=&\textcolor{blue}{\sum_{\substack{j\neq m}}^{M}\left[\Phi_\kappa(z_m,z_j) Q_j |\partial D_m|+O\left(\frac{a^{5-\beta}}{d}\left(\kappa+\frac{1}{d}\right) \right)\right].}\nonumber\\
 &=&\sum_{\substack{j\neq m}}^{M}\left[\Phi_\kappa(z_m,z_j) Q_j |\partial D_m|+O\left(\frac{a}{d_{mj}}\left(\kappa+\frac{1}{d_{mj}}\right)|\partial D_m| \left|\int_{\partial D_j} \sigma_{j} (s)ds\right|\right)\right]\nonumber\\
 &=&\sum_{\substack{j\neq m}}^{M}\left[\Phi_\kappa(z_m,z_j) Q_j |\partial D_m|+O\left(\frac{a^{5-\beta}}{d_{mj}}\left(\kappa+\frac{1}{d_{mj}}\right) \right)\right]\nonumber\\
  &=&\sum_{\substack{j\neq m}}^{M}\Phi_\kappa(z_m,z_j) Q_j |\partial D_m|+\sum_{n=1}^{[d^{-\alpha}]}[(2n+1)^3-(2n-1)^3]O\left(\frac{a^{5-\beta}}{nd^{\alpha}}\left(\kappa+\frac{1}{ nd^{\alpha}}\right) \right)\nonumber\\
  &=&\sum_{\substack{j\neq m}}^{M}\Phi_\kappa(z_m,z_j) Q_j |\partial D_m|+O\left(2\frac{a^{5-\beta}}{d^{\alpha}}\sum_{n=1}^{[d^{-\alpha}]}\frac{12n^2+1}{n}\left(\kappa+\frac{1}{ nd^{\alpha}}\right) \right)\nonumber\\
  &=&\sum_{\substack{j\neq m}}^{M}\Phi_\kappa(z_m,z_j) Q_j |\partial D_m|+O\left(2\frac{a^{5-\beta}}{d^{2\alpha}}\left[7\kappa+\frac{6\kappa+13}{d^{\alpha}} \right] \right)
  \end{eqnarray}
  and
\begin{eqnarray}\label{approximating-Dpr}
|D{^\prime}|&=&\left\vert\sum_{j\neq m}^{M}\int_{\partial D_j}\left(\int_{\partial D_m}\left[\Phi_\kappa(s_m,s)-\Phi_\kappa(s_m,z_j)\right] d{s_m}\right) \sigma_{j} (s)ds\right\vert\nonumber\\
% % %  &=&\textcolor{blue}{\sum_{j\neq m}^{M}O\left(\frac{a^{5-\beta}}{d}\left(\kappa+\frac{1}{d}\right)\right)}\nonumber\\
&=&\sum_{j\neq m}^{M}O\left(\frac{a^{5-\beta}}{d_{mj}}\left(\kappa+\frac{1}{d_{mj}}\right)\right)\nonumber\\
&=&\sum_{n=1}^{[d^{-\alpha}]}[(2n+1)^3-(2n-1)^3]O\left(\frac{a^{5-\beta}}{nd^{\alpha}}\left(\kappa+\frac{1}{nd^{\alpha}}\right)\right)\nonumber\\
&=&O\left(2\frac{a^{5-\beta}}{d^{\alpha}}\sum_{n=1}^{[d^{-\alpha}]}\frac{12n^2+1}{n}\left(\kappa+\frac{1}{nd^{\alpha}}\right)\right)\nonumber\\
&=&O\left(2\frac{a^{5-\beta}}{d^{2\alpha}}\left[7\kappa+\frac{6\kappa+13}{d^{\alpha}} \right]\right).
\end{eqnarray}
\item Now, let us approximate $C$. Since $\vert \int_{\partial D_m}\Phi_{\kappa}(s_m, s)ds_m\vert \leq \frac{1}{4\pi}\int_{\partial D_m}\frac{1}{\vert s_m- s\vert}ds_m =O(a)$, then

\begin{eqnarray}\label{approximating-C1}
|C|&=&\left\vert\int_{\partial D_m}\left(\int_{\partial D_m}\Phi_\kappa(s_m,s) d{s_m}\right) \sigma_{m} (s)ds\right\vert\nonumber\\
&=&O\left( a^{3-\beta} \right).
\end{eqnarray}

\end{itemize}

Hence, from (\ref{Q_mintdbl-1}-\ref{approximating-C1}), we obtain the approximation below

\begin{eqnarray}\label{Q_mintdbl-3}
-Q_m&+&\sum_{\substack{j\neq m}}^{M}\Phi_\kappa(z_m,z_j)\lambda_m|\partial D_m|Q_j\nonumber\\
&=&-\lambda_m |\partial D_m| e^{i\kappa\theta\cdot{z_m}}+O\left((\vert\lambda_m\vert+\kappa)\kappa a^3\right)+\lambda_m O\left(a^{3-\beta}\right),\nonumber\\
&&\quad+O\left(\kappa^2 a^{4-\beta}\right)+O\left(2\kappa^2\frac{a^{5-\beta}}{d^{2\alpha}}[\frac{6}{d^{\alpha}}+7]\right)+\lambda_m O\left(2\frac{a^{5-\beta}}{d^{2\alpha}}\left[7\kappa+\frac{6\kappa+13}{d^{\alpha}} \right]\right).% \nonumber\\
\end{eqnarray}
Indeed, $\int_{\partial D_m}\lambda_m ((U^i{(s_m)}-U^{i}(z_m))d{s_m}=O(\vert\lambda\vert\kappa a^3)$ and $\int_{\partial D_m}\frac{\partial U^i}{\partial \nu_m}{(s_m)} d{s_m}=O(\kappa^2 a^3)$. %Here $J_m$ is the constant defined as in \eqref{defofJm}. 
In addition, since $\beta <1$, then  $\lambda_m O\left(a^{3-\beta}\right)=O(a^{3-2\beta})=o(a^{2-\beta})=o(\lambda |\partial D_m| e^{i\kappa\theta\cdot{z_m}})$.

\par 
We can rewrite the above algebraic system as
\begin{eqnarray}\label{Q_mintdbl-4f}
-\frac{1}{\lambda_m|\partial D_m|} Q_m
&=&- e^{i\kappa\theta\cdot{z_m}}-\sum_{\substack{j\neq m}}^{M}\Phi_\kappa(z_m,z_j)Q_j+ Err
\end{eqnarray}
with 
\begin{eqnarray}\label{Err}
 Err
&:=&\frac{1}{\lambda_m|\partial D_m|}\left[O\left((\lambda_+ +\kappa\epsilon^\beta)\kappa a^{3-\beta}\right)+O\left(\kappa^2 a^{4-\beta}\right)+O\left(2\kappa^2\frac{a^{5-\beta}}{d^{2\alpha}}[\frac{6}{d^{\alpha}}+7]\right)\right.\nonumber\\
&&\left.+ O\left(2\frac{a^{5-2\beta}}{d^{2\alpha}}\left[7\kappa+\frac{6\kappa+13}{d^{\alpha}} \right]\right)+O\left(a^{3-2\beta}\right)\right] \nonumber\\
&=&\frac{1}{\lambda_m|\partial D_m|}O\left( a^{3-2\beta} +\frac{a^{5-2\beta}}{d^{3\alpha}}\right)\nonumber\\
 &=&O\left( a^{1-\beta}+\frac{a^{3-\beta}}{d^{3\alpha}}\right).
 \end{eqnarray}
In the last two lines in the above approximation, we used the fact that $\kappa\leq\kappa_{\max}$ and $d\leq d_{\max}$.

\end{proofn}
%\subsubsection{Invertibility of the the algebraic system}
\subsubsection{Invertibility of the algebraic system}\label{algsystem-sec-mainthrmproof-small}
We define the following algebraic system

\begin{eqnarray}\label{Q_mintdbl-algsys}
\frac{\bar{Q}_m}{\bar{C}_m} 
&:=&- U^i(z_m) -\sum_{\substack{j\neq m}}^{M}\bar{C}_j \Phi_\kappa(z_m,z_j)\frac{\bar{Q}_j}{\bar{C}_j}
\end{eqnarray}
 for all $m=1,2,\dots,M$. %It is valid with the error of order $O\left(\kappa a+\frac{\kappa a^2}{d}+\frac{a^2}{d^2}\right)$.
It can be written in a compact form as
\begin{equation}\label{compacfrm1}
 \mathbf{B}\bar{Q}=\mathrm{U}^I,
\end{equation}
\noindent
where $\bar{Q},\mathrm{U}^I \in \mathbb{C}^{M\times 1}\mbox{ and } \mathbf{B}\in\mathbb{C}^{M\times M}$ are defined as
\begin{eqnarray}
\mathbf{B}:=\left(\begin{array}{ccccc}
   -\frac{1}{\bar{C}_1} &-\Phi_\kappa(z_1,z_2)&-\Phi_\kappa(z_1,z_3)&\cdots&-\Phi_\kappa(z_1,z_M)\\
-\Phi_\kappa(z_2,z_1)&-\frac{1}{\bar{C}_2}&-\Phi_\kappa(z_2,z_3)&\cdots&-\Phi_\kappa(z_2,z_M)\\
 \cdots&\cdots&\cdots&\cdots&\cdots\\
-\Phi_\kappa(z_M,z_1)&-\Phi_\kappa(z_M,z_2)&\cdots&-\Phi_\kappa(z_M,z_{M-1}) &-\frac{1}{\bar{C}_M}
   \end{array}\right),\label{mainmatrix-acoustic-small}\\
\nonumber\\
%\begin{split}
 \bar{Q}:=\left(\begin{array}{cccc}
    \bar{Q}_1 & \bar{Q}_2 & \ldots  & \bar{Q}_M
   \end{array}\right)^\top \text{ and } 
\mathrm{U}^I:=\left(\begin{array}{cccc}
     U^i(z_1) & U^i(z_2)& \ldots &  U^i(z_M)
   \end{array}\right)^\top.\label{coefficient-and-incidentvectors-acoustic-small}
%\end{split}
\end{eqnarray}
The above linear algebraic system is solvable for $\bar{Q}_j,~1\leq j\leq M$, when the matrix $\mathbf{B}$ is invertible.
% To use the simple notations, let us denote $\frac{\lambda|\partial D_m|}{\left[-1+\lambda J_m\right]}$ by $\bar{C}_{m}$ for $m=1,2,\dots,M$. 
Now we give sufficient conditions for the invertibility of system \eqref{compacfrm1}. We recall that $\lambda_m:=\lambda_{m,0}a^{-\beta},\; \beta \geq 0$.

\begin{lemma}\label{Mazyawrkthm} We distinguish the following two cases:

\begin{itemize}
 \item Let $\Re (\lambda_{m,0}) < 0$ \footnote{  In this case, we can actually relax the condition $s \leq 2- \beta$ on the number of small scatterers $s$, see Remark \ref{Remark-on-sign-lambda}} and assume that 
 $\min\limits_{1\leq{m}\leq{M}}\frac{\Re\bar{C}_m}{\vert\bar{C}_m\vert^2}>\frac{\sqrt{2 M_{max}}}{\pi\,d^{\frac{s}{t}}}$ 
 then the matrix $\mathbf{B}$ 
is invertible and the solution vector $\bar{Q}$ of \eqref{compacfrm1} satisfies the estimate
\begin{equation}\label{mazya-fnlinvert-small-ac-2}
\begin{split}
 \sum_{m=1}^{M}|\bar{Q}_m|^{2}
\leq4 \left(\frac{\min\limits_{1\leq{m}\leq{M}}\Re\bar{C}_m}{\max\limits_{1\leq{m}\leq{M}}\vert\bar{C}_m\vert^2}-\frac{\sqrt{2 M_{max}}}{\pi\,d^{\frac{s}{t}}}\right)^{-2}\sum_{m=1}^{M}\left|U^i(z_m)\right|^2.
\end{split}
\end{equation}

\item Let $\Re (\lambda_{m,0}) > 0$ and assume that $\frac{\min\limits_{1\leq{m}\leq{M}}\Re(-\bar{C}_m)}{(\max\limits_{1\leq{m}\leq{M}}\vert\bar{C}_m\vert)^2}>\frac{\sqrt{2M_{max}}}{\pi d^{\frac{s}{t}}}$ then
 the matrix $\mathbf{B}$ 
is invertible and the solution vector $\bar{Q}$ of \eqref{compacfrm1} satisfies the estimate
\begin{equation}\label{mazya-fnlinvert-small-ac-2-1}
 \begin{split}
 \sum_{m=1}^{M}|\bar{Q}_m|^{2}
\leq 4\left(\frac{\min\limits_{1\leq{m}\leq{M}}\Re(-\bar{C}_m)}{\max\limits_{1\leq{m}\leq{M}}\vert\bar{C}_m\vert^2}-\frac{\sqrt{2 M_{max}}}{\pi\,d^{\frac{s}{t}}}\right)^{-2}\sum_{m=1}^{M}\left|U^i(z_m)\right|^2.
 \end{split}
 \end{equation}
 where $M_{max}:=M\; a^s$, recalling that $M:=O(a^{-s})$, as $a\rightarrow 0$.
 \end{itemize}
\end{lemma}
Since $\bar{C}_m:=-\lambda_m \vert \partial D_m \vert$ and $d\geq d_{min} a^t$, 
the condition $\frac{\min\limits_{1\leq{m}\leq{M}}\vert\Re\bar{C}_m\vert}{\max\limits_{1\leq{m}\leq{M}}\vert\bar{C}_m\vert^2}>\frac{\sqrt{2M_{max}}}{\pi\,d^{\frac{s}{t}}}$ 
is satisfied if $\frac{4\lambda_- a^{s-(2-\beta)}}{\pi \lambda^2_+}>\frac{\sqrt{2 M_{max}}}{\pi d^{\frac{s}{t}}}$. This is possible if $s<2-\beta$ and $a$ is small enough or
$s=2-\beta$ and $\lambda_-$ and $\lambda_+$ satisfy $\frac{4\lambda_-}{\pi \lambda^2_+}>\frac{\sqrt{2 M_{max}}}{\pi d_{min}^{\frac{s}{t}}}$. Recall that $t\geq \frac{s}{3}$, hence 
$\frac{s}{t}\leq 3$. If we assume $d_{min}\leq 1$, then the last condition is satisfied if $\frac{4\lambda_-}{\pi \lambda^2_+}>\frac{\sqrt{2 M_{max}}}{\pi d_{min}^{3}}$ and if
 $d_{min}\geq 1$ then we take $\frac{4\lambda_-}{\pi \lambda^2_+}>\frac{\sqrt{2 M_{max}}}{\pi}$.
The proof of this lemma will be given in the appendix. 

%\begin{proofn}[\textbf{End of the Proof of Theorem \ref{Maintheorem-ac-small-sing} :} ]\\
%\ ~ \ \par
\begin{proofe}
We can rewrite the inequatilies \eqref{mazya-fnlinvert-small-ac-2} and \eqref{mazya-fnlinvert-small-ac-2-1} using norm inequalities as
\begin{eqnarray}\label{mazya-fnlinvert-small-ac-3}
 \sum_{m=1}^{M}|\bar{Q}_m|
&\leq&2\left(\frac{\min\limits_{1\leq{m}\leq{M}}\Re\bar{C}_m}{\max\limits_{1\leq{m}\leq{M}}\vert\bar{C}_m\vert}-\frac{\sqrt{2M_{max}}}{\pi\,d^{\frac{s}{t}}}\max\limits_{1\leq m \leq M}{\vert\bar{C}_m\vert}\right)^{-1}M\max\limits_{1\leq m \leq M}|\bar{C}_m|\max\limits_{1\leq m \leq M}\left|U^i(z_m)\right|
\\\label{mazya-fnlinvert-small-ac-3-1}
\sum_{m=1}^{M}|\bar{Q}_m|
&\leq&\left(\frac{\min\limits_{1\leq{m}\leq{M}}\Re (-\bar{C}_m)}{\max\limits_{1\leq{m}\leq{M}}\vert\bar{C}_m\vert}-\frac{\sqrt{2M_{max}}}{\pi\,d^{\frac{s}{t}}}\max\limits_{1\leq m \leq M}{\vert\bar{C}_m\vert}\right )^{-1}M\max\limits_{1\leq m \leq M}|\bar{C}_m|\max\limits_{1\leq m \leq M}\left|U^i(z_m)\right|
\end{eqnarray}
which holds for the cases $\Re \lambda_m \leq 0$ and $\Re \lambda_m \geq 0$ respectively.
\bigskip

\subsubsection{The dominant part of the total charges}
The difference between \eqref{Q_mintdbl-algsys} and \eqref{Q_mintdbl-4} produce the following
  \begin{eqnarray}\label{qcdiftilde}
   \frac{{Q}_m-\bar{Q}_m}{\bar{C}_m} &=&-\sum_{\substack{j=1 \\ j\neq m}}^{M} \Phi_\kappa(z_m,z_j)\left({Q}_j-\bar{Q}_j\right)+Err,
  \end{eqnarray}
for $m=1,\dots,M$. Comparing the above system of equations \eqref{qcdiftilde} with \eqref{Q_mintdbl-algsys} and by making use of the estimates \eqref{mazya-fnlinvert-small-ac-3} and \eqref{mazya-fnlinvert-small-ac-3-1}, we obtain
 \begin{eqnarray}\label{unncmaybe}
\sum_{m=1}^{M}({{Q}_m-\bar{Q}_m})&=& O\left(M{a^{2-\beta}}Err\right).
\end{eqnarray}
We can evaluate the $\bar{Q}_m$'s from the algebraic system 
\eqref{Q_mintdbl-algsys}. This means that $\bar{Q}_m$'s are the dominant parts 
of the total charges $Q_m$'s.

\subsection{Approximation of the far-fields. III. End of the proof of Theorem \ref{Maintheorem-ac-small-sing}}\label{sec-mainthrmproof-small} Using  \eqref{unncmaybe} in 
\eqref{x oustdie1 D_m} we can represent the far-field pattern in terms of 
$\bar{Q}_m$ as follows
\begin{eqnarray}\label{x oustdie1 D_m1}
\hspace{-.7cm}U^{\infty}(\hat{x})&=&\sum_{m=1}^{M}[e^{-i\kappa\hat{x}\cdot z_{m}}Q_m+O(\kappa\,a^{3-\beta})]\nonumber\\
                   &=&\sum_{m=1}^{M}\left[e^{-i\kappa\hat{x}\cdot z_{m}}[\bar{Q}_m+({Q}_m-\bar{Q}_m)]+O(\kappa\,a^{3-\beta})\right]\nonumber\\
 &=&\sum_{m=1}^{M}e^{-i\kappa\hat{x}\cdot z_{m}}\bar{Q}_m+O(M{a^{2-\beta}}[Err+(\kappa\,a)])\nonumber\\
 &\substack{=\\ \eqref{Err}}&\sum_{m=1}^{M}e^{-i\kappa\hat{x}\cdot z_{m}}\bar{Q}_m+O\left(M{a^{2-\beta}}\left(  a^{1-\beta}%+\frac{a^{3}}{d^{4\alpha}}
+\frac{a^{3-\beta}}{d^{3\alpha}}\right)\right)\nonumber\\
&=& \sum_{m=1}^{M}e^{-i\kappa\hat{x}\cdot z_{m}}\bar{Q}_m+O\left(a^{3-s-2\beta}\right)
\end{eqnarray}
since, as $s=3 t \alpha \leq 2-\beta$, we have $\frac{a^{3-\beta}}{d^{3\alpha}}\sim a^{3-\beta-3t\alpha}=O(a)=\circ(a^{1-\beta})$ if $\beta <1$. Hence Theorem \ref{Maintheorem-ac-small-sing} is proved with the replacement of ${C}_m$ in the statement by $\bar{C}_m$. 

% \footnote{If we do not distinguish the near by and far obstacles, as it discussed immediately after \eqref{nrm1singulayer2dd-ccc}, and estimating the terms $B, D, B^{\prime}$ and $D^{\prime}$, % \cite{DPC-SM:MMS2013,DPC-SM:MANA2015} 
% we observe the behaviour of $Err$ term as   $Err=O\left( a^{1-\beta} +M\frac{a^{3-\beta}}{d^2}\right)$. Then the approximation error of the farfield expression \eqref{x oustdie1 D_m1} should be replaced by %Theorem \ref{Maintheorem-ac-small-sing} 
% $O\left(Ma^{2-\beta}\left( a^{1-\beta} +M\frac{a^{3-\beta}}{d^2}\right)\right)$. }
\end{proofe}

\section{Proof of Corollary \ref{Maintheorem-ac-small-sing-1}}
The steps of the proof of Corollary \ref{Maintheorem-ac-small-sing-1} are the same as for the proof of Theorem \ref{Maintheorem-ac-small-sing}. Here we only explain the main changes that are needed to derive it. 

We recall that, for an obstacle $D_{\epsilon}$ of radius $\epsilon$, $\mathcal{S}(\phi)(s):=\int_{\partial D_{\epsilon}}\Phi_{\kappa}(s, t)\phi(t)dt$ and $\mathcal{D}(\phi)(s):=\int_{\partial D_{\epsilon}}\frac{\partial \Phi_{\kappa}(s, t)}{\partial \nu(t)}\phi(t)dt$. Similarly, we set $\mathcal{S}_G(\phi)(s):=\int_{\partial D_{\epsilon}}G_{\kappa}(s, t)\phi(t)dt$ and $\mathcal{D}_G(\phi)(s):=\int_{\partial D_{\epsilon}}\frac{\partial G_{\kappa}(s, t)}{\partial \nu(t)}\phi(t)dt$.

We see that $W_{\kappa}(x, z):=G_{\kappa}(x, z)-\Phi_{\kappa}(x, z)$ satisfies
\begin{equation}
(\Delta +\kappa^2n^2)W_{\kappa}=\kappa^2 (1-n^2) \Phi_{\kappa},\; \mbox{ in } \mathbb{R}^3
\end{equation}
with the Sommerfeld radiation conditions.  Since $\Phi_{\kappa}(\cdot, z)$, $z \in \mathbb{R}^3$ is bounded in $L^p(\Omega)$, for $p<3$, by interior estimates, we deduce that $W(\cdot, z)$, $z \in \mathbb{R}^3$ is bounded in $W^{2, p}(\Omega)$, for $p<3$, and hence, in particular, the normal traces are bounded in $L^2(\partial D_{\epsilon})$. Then we can show that the norms of the operators
\begin{equation}\label{single-layer-G_k}
\mathcal{S}_G- \mathcal{S}: L^2(\partial D_{\epsilon}) \rightarrow H^1(\partial D_{\epsilon})
\end{equation}   
and 
\begin{equation}\label{double-layer-G_k}
\mathcal{D}_G- \mathcal{D}: L^2(\partial D_{\epsilon}) \rightarrow L^2(\partial D_{\epsilon})
\end{equation}
are of the order $O(\epsilon)$ at least. The representation (\ref{qcimprequiredfrm1}), in Proposition \ref{existence-of-sigmas}, needs to be replaced by
 \begin{equation}\label{qcimprequiredfrm1-G_k}
  U_n^{t}(x)=V_n^{t}(x)+\sum_{m=1}^{M}\int_{\partial D_m}G_\kappa(x,s)\sigma_{m} (s)ds,~x\in\mathbb{R}^{3}\backslash\left(\mathop{\cup}_{m=1}^M \bar{D}_m\right), 
\end{equation}
where $V_n^{t}$ is the total field corresponding to the background modeled by the index of refraction $n$, see (\ref{acimpoenetrable-2}). We use the single layer potentials defined by the Greens' function $G_{\kappa}$ instead of the fundamental function $\Phi_{\kappa}$. The main tools used in justifying Proposition \ref{existence-of-sigmas} are the invertibility properties of the corresponding integral operators (i.e. the Fredholm property) and the jumps of the double layer potentials defined by $\Phi_{\kappa}$. These two tools are satisfied also when we use $G_{\kappa}$ instead of $\Phi_{\kappa}$ due to the error bounds in (\ref{single-layer-G_k}) and (\ref{double-layer-G_k}).

The a priori estimates on the densities $\sigma_m$'s derived in section \ref{SLPR-2} are  quantitative versions of the result in Proposition \ref{existence-of-sigmas}. Due to  error bounds in (\ref{single-layer-G_k}) and (\ref{double-layer-G_k}), those estimates can then be translated to the densities used in the representation (\ref{qcimprequiredfrm1-G_k}). Of course, in Lemma \ref{rep1singulayer} one needs to replace the equalities by inequalities to estimate the properties of the operators defined by $G_{\kappa}$, on the scaled obstacles $D_{\epsilon}$, in terms of the properties of the operators defined by $\Phi_{\kappa}$ on the original obstacles $B$.

Finally, to do the same analysis as in section \ref{SLPR-3}, one needs to split $G_{\kappa}$ as $G_{\kappa}=\Phi_{\kappa}+(G_{\kappa}-\Phi_{\kappa})$ and use the results derived in section \ref{SLPR-3} and again error bounds in (\ref{single-layer-G_k}) and (\ref{double-layer-G_k}).

\section{Justification of Remark \ref{Remark-spherical-scatterers}}
We show only the main changes needed in the proof of Theorem \ref{Maintheorem-ac-small-sing}. This change occurs
in the proof of Proposition \ref{fracqfracc-ac} and precisely in evaluating the term $C$, i.e. (\ref{approximating-C1}). We rewrite $C$ as
\begin{equation}\label{approximating-C1-sphere}
C = \int_{\partial D_m}\left(\int_{\partial D_m}\Phi_0(s_m,s) d{s_m}\right) \sigma_{m} (s)ds+  \int_{\partial D_m}\left(\int_{\partial D_m}\Phi_\kappa(s_m,s)-\Phi_0(s_m,s)d{s_m}\right) \sigma_{m} (s)ds. \nonumber
\end{equation} 
Since 
\begin{eqnarray*}
\left\vert \int_{\partial D_m}\left(\int_{\partial D_m}\left[\Phi_\kappa(s_m,s)-\Phi_0(s_m,s)\right]) d{s_m}\right) \sigma_{m} (s)ds \right\vert
&<&\frac{1}{2}\sum_{l=1}^\infty \frac{\kappa^la^{l-1}}{l!} \vert\partial D_m\vert  \int_{\partial D_m}\vert \sigma_{m} (s)\vert ds\\
&=&O\left( \kappa a^{4-\beta} \right)
\end{eqnarray*}
then we can write $C$ as follows;
\begin{eqnarray}\label{approximating-C}
C&=&\int_{\partial D_m}\left(\int_{\partial D_m}\Phi_0(s_m,s) d{s_m}\right) \sigma_{m} (s)ds +O\left( \kappa a^{4-\beta} \right)\nonumber\\
&=&I_mQ_m + O\left( \kappa a^{4-\beta} \right),
\end{eqnarray}
where $I_m=\int_{\partial D_m}\Phi_0(s_m,t) d{s_m},$ is a constant for each $t\in\partial D_m$ if $D_m$ is a ball.  Indeed, for the single layer potential $\left(\mathcal{S}^{i_{\kappa}}_{ D_m}\psi\right) (s)=\int_{\partial D_\epsilon}\Phi_0(s,t)\psi(t)dt$, we have the jump condition
as $
 \left.\frac{\partial }{\partial \nu}\left(\mathcal{S}^{i_{\kappa}}_{ D_m}\psi\right) (s)\right|_{-}=\left(\frac{1}{2}\textbf{I}+\mathbf{K}^{i_{\kappa}^{*}}_{ D_m}\right)\psi (s)
$
recalling that $\mathbf{K}^{i_{\kappa}^{*}}_{ D_m}$ is the adjoint of the double layer operator. Since $\mathbf{K}^{i_{\kappa}^{*}}_{ D_m}[1]=\mathbf{K}^{i_{\kappa}}_{ D_m}[1]=-1/2$, then
we have $\frac{\partial }{\partial \nu^{-}}\mathcal{S}^{i_{\kappa}}_{D_m} [1]=0$ on $\partial D_m$.  
Hence $\mathcal{S}^{i_{\kappa}}_{D_m}[1]$ is a constant, namely $I_m$, as it satisfies $\Delta \;\mathcal{S}^{i_{\kappa}}_{D_m}[1]=0\; \mbox{ in } D_m$ and has zero normal derivative on $\partial D_m$.

\bigskip

With this correction at hand the estimate (\ref{Q_mintdbl-3}) becomes

\begin{eqnarray}\label{Q_mintdbl-3-sphere}
(-1+\lambda_m I_m)Q_m&+&\sum_{\substack{j\neq m}}^{M}\Phi_\kappa(z_m,z_j)\lambda_m|\partial D_m|Q_j\nonumber\\
&=&-\lambda_m |\partial D_m| e^{i\kappa\theta\cdot{z_m}}+O\left((\vert\lambda_m\vert+\kappa)\kappa a^3\right)+\lambda_m O\left(a^{4-\beta}\right),\nonumber\\
&&\quad+O\left(\kappa^2 a^{4-\beta}\right)+O\left(2\kappa^2\frac{a^{5-\beta}}{d^{2\alpha}}[\frac{6}{d^{\alpha}}+7]\right)+\lambda_m O\left(2\frac{a^{5-\beta}}{d^{2\alpha}}\left[7\kappa+\frac{6\kappa+13}{d^{\alpha}} \right]\right).% \nonumber\\
\end{eqnarray}

and the system (\ref{Q_mintdbl-4f}) as
\begin{eqnarray}\label{Q_mintdbl-4f-sphere}
\frac{-1+\lambda_m I_m}{\lambda_m|\partial D_m|} Q_m
&=&- e^{i\kappa\theta\cdot{z_m}}-\sum_{\substack{j\neq m}}^{M}\Phi_\kappa(z_m,z_j)Q_j+ Err
\end{eqnarray}
with 
\begin{eqnarray}\label{Err1}
 Err
&:=&\frac{1}{\lambda_m|\partial D_m|}\left[O\left((\lambda_+ +\kappa\epsilon^\beta)\kappa a^{3-\beta}\right)+O\left(\kappa^2 a^{4-\beta}\right)+O\left(2\kappa^2\frac{a^{5-\beta}}{d^{2\alpha}}[\frac{6}{d^{\alpha}}+7]\right)\right.\nonumber\\
&&\left.+ O\left(2\frac{a^{5-2\beta}}{d^{2\alpha}}\left[7\kappa+\frac{6\kappa+13}{d^{\alpha}} \right]\right)+O\left(a^{4-2\beta}\right)\right] \nonumber\\
&=&\frac{1}{\lambda_m|\partial D_m|}O\left( a^{3-\beta} +\frac{a^{5-2\beta}}{d^{3\alpha}}\right)\nonumber\\
 &=&O\left( a+\frac{a^{3-\beta}}{d^{3\alpha}}\right).
 \end{eqnarray}

Finally, the estimate (\ref{x oustdie1 D_m1}) becomes

\begin{eqnarray}\label{x oustdie1 D_m1-sphere}
\hspace{-.7cm}U^{\infty}(\hat{x})&=&\sum_{m=1}^{M}[e^{-i\kappa\hat{x}\cdot z_{m}}Q_m+O(\kappa\,a^{3-\beta})]\nonumber\\
                   &=&\sum_{m=1}^{M}\left[e^{-i\kappa\hat{x}\cdot z_{m}}[\bar{Q}_m+({Q}_m-\bar{Q}_m)]+O(\kappa\,a^{3-\beta})\right]\nonumber\\
 &=&\sum_{m=1}^{M}e^{-i\kappa\hat{x}\cdot z_{m}}\bar{Q}_m+O(M{a^{2-\beta}}[Err+(\kappa\,a)])\nonumber\\
 &\substack{=\\ \eqref{Err1}}&\sum_{m=1}^{M}e^{-i\kappa\hat{x}\cdot z_{m}}\bar{Q}_m+O\left(M{a^{2-\beta}}\left(  a%+\frac{a^{3}}{d^{4\alpha}}
+\frac{a^{3-\beta}}{d^{3\alpha}}\right)\right)\nonumber\\
&=& \sum_{m=1}^{M}e^{-i\kappa\hat{x}\cdot z_{m}}\bar{Q}_m+O\left(a^{3-s-\beta}\right)
\end{eqnarray}
since, as $s=3 t \alpha \leq 2-\beta$, we have $\frac{a^{3-\beta}}{d^{3\alpha}}=O(a^{3-\beta-3t\alpha})=O(a)$.% and $a^{2-\beta}=O(a)$ if $\beta \leq 1$.

 \section{Appendix: \it{Proof of Lemma \ref{Mazyawrkthm}}}
\begin{proofe}

We start by factorizing $\mathbf{B}$ as $\mathbf{B}=-(\mathbf{C}^{-1}+\mathbf{B}_{n})$ where $\mathbf{C}:=Diag(\bar{C}_1,\bar{C}_2,\dots,\bar{C}_M)\in\mathbb{R}^{M\times\,M}$, $I$ is the identity matrix and 
$\mathbf{B}_{n}:=-\mathbf{C}^{-1}-\mathbf{B}$. We have $\mathbf{B}:\mathbb{C}^{M}\rightarrow\mathbb{C}^{M}$, so it is enough to prove the injectivity in order to prove its invetibility. For this purpose,  let $X,Y$ are vectors in $\mathbb{C}^{M}$ and
 consider the system
\begin{eqnarray}\label{systemsolve1-small-ac}(\mathbf{C}^{-1}+\mathbf{B}_{n})X&=&Y.\end{eqnarray}
Let ${(\cdot)}^{real}$ and ${(\cdot)}^{img}$ denotes the real and the imaginary parts of the corresponding complex number/vector/matrix. For convenience, let us denote $\mathbf{C}^{-1}$ by $\mathbf{C}_{I}$. Now, the following can be written from \eqref{systemsolve1-small-ac};
\begin{eqnarray}
 (\mathbf{C}_{I}^{real}+\mathbf{B}^{real}_{n})X^{real}-(\mathbf{C}_{I}^{img}+\mathbf{B}^{img}_{n})X^{img}&=&Y^{real}\label{systemsolve1-small-sub1-ac-1},\\
 (\mathbf{C}_{I}^{real}+\mathbf{B}^{real}_{n})X^{img}+(\mathbf{C}_{I}^{img}+\mathbf{B}^{img}_{n})X^{real}&=&Y^{img}\label{systemsolve1-small-sub2-ac-1},
\end{eqnarray}
which leads to
\begin{eqnarray}
 \langle\,(\mathbf{C}_{I}^{real}+\mathbf{B}^{real}_{n})X^{real},X^{real}\rangle\,-\langle\,(\mathbf{C}_{I}^{img}+\mathbf{B}^{img}_{n})X^{img},X^{real}\rangle&=&\langle\,Y^{real},X^{real}\rangle\label{systemsolve1-small-sub1-ac-2},\\
 \langle\,(\mathbf{C}_{I}^{real}+\mathbf{B}^{real}_{n})X^{img},X^{img}\rangle\,+\langle\,(\mathbf{C}_{I}^{img}+\mathbf{B}^{img}_{n})X^{real},X^{img}\rangle&=&\langle\,Y^{img},X^{img}\rangle\label{systemsolve1-small-sub2-ac-2}.
\end{eqnarray}
By summing up \eqref{systemsolve1-small-sub1-ac-2} and \eqref{systemsolve1-small-sub2-ac-2} will give
\begin{equation}
\begin{split}
\langle\,\mathbf{C}_{I}^{real}X^{real},X^{real}\rangle\,+\langle\,\mathbf{B}^{real}_{n}X^{real},X^{real}\rangle\,+\langle\,\mathbf{C}_{I}^{real}X^{img},X^{img}\rangle\,+\langle\,\mathbf{B}^{real}_{n}X^{img},X^{img}\rangle\,\\
% \langle\,(I+\mathbf{B}^{real}_{n}\mathbf{C})X^{real},\mathbf{C}X^{real}\rangle\,+\langle\,(I+\mathbf{B}^{real}_{n}\mathbf{C})X^{img},\mathbf{C}X^{img}\rangle\,
=\langle\,Y^{real},X^{real}\rangle+\langle\,Y^{img},X^{img}\rangle.\label{systemsolve1-small-sub1-2-ac-3}
\end{split}
\end{equation}
We can observe that, the right-hand side in \eqref{systemsolve1-small-sub1-2-ac-3} does not exceed
\begin{equation}\label{systemsolve1-small-sub1-2-ac-4}
\begin{split}
\langle\,X^{real},X^{real}\rangle^{1\slash\,2}\langle\,Y^{real},Y^{real}\rangle^{1\slash\,2}+\langle\,X^{img},X^{img}\rangle^{1\slash\,2}\langle\,Y^{img},Y^{img}\rangle^{1\slash\,2}\\
\leq2\langle\,X^{|\cdot|},X^{|\cdot|}\rangle^{1\slash\,2}\langle\,Y^{|\cdot|},Y^{|\cdot|}\rangle^{1\slash\,2}.
\end{split}
\end{equation}
At this stage, we divide the proof into two cases. 

\begin{enumerate}
 \item If $\Re (\lambda_{m, 0}) < 0$ for each $m$. In this case $\Re\bar{C}_m > 0$.
 We know that:
\begin{equation}
\vert \langle\; B_n^{real} X^{real},\; X^{real}\rangle \vert \leq \Vert B_n^{real}\Vert_2 \vert X^{real}\vert^2_2
\end{equation}
where $\Vert B_n^{real}\Vert^2_2:=\sum^M_{i,\; j =1}(B^{real}_n)^2_{i,j}$ and $(B^{real}_n)_{i,j}:=\Re\; \Phi(z_i, z_j)$ if $i\ne j$ and $(B^{real}_n)_{i,i}:=0$ for $i,j=1,...,M$. Hence $\vert (B^{real}_n)_{i,j} \vert \leq \frac{1}{4\pi \vert z_i-z_j\vert },\; i\ne j$. Arguing as in (\ref{KDKnrm1}), we see that
\begin{equation}
\sum^M_{i, j=1}(B^{real}_n)^2_{i,j}\leq M\sum^{[d^{-\alpha}]}_{n=1}[(2n+1)^3-(2n-1)^3]\frac{1}{(4\pi)^2 n^2 d^{2\alpha}}\leq M \frac{2d^{-2\alpha}}{\pi^2}\sum^{[d^{-\alpha}]}_{n=1}1=\frac{2M d^{-3\alpha}}{\pi^2}.
\end{equation}
Observing that $M=O(a^{-s})\leq M_{max}a^{-s}$ and that $s=3t\alpha$, we obtain:
\begin{equation}\label{estimate-B-n-negative-realpart}
\Vert B_n^{real}\Vert_2 \leq \frac{\sqrt{2M_{max}}}{\pi}a^{-\frac{s}{t}}.
\end{equation}
From (\ref{systemsolve1-small-sub1-2-ac-3}) and (\ref{systemsolve1-small-sub1-2-ac-4}), we deduce that
 \begin{equation}\label{fnlinvert-small-ac-2-f1-1}
\begin{split}
 \left(\frac{\min\limits_{1\leq{m}\leq{M}}\Re\bar{C}_m}{(\max\limits_{1\leq{m}\leq{M}}\vert\bar{C}_m\vert)^2}-\frac{\sqrt{2M_{max}}}{\pi\,d^{\frac{s}{t}}}\right)\sum_{m=1}^{M}|X_m|^{2}
\leq2\left(\sum_{m=1}^{M}|X_m|^2\right)^{1/2}\left(\sum_{m=1}^{M}|Y_m|^2\right)^{1/2},
\end{split}
\end{equation}
which yields
\begin{equation}\label{fnlinvert-small-ac-2-f1-2}
\begin{split}
 \sum_{m=1}^{M}|X_m|^{2}
\leq 4\left(\frac{\min\limits_{1\leq{m}\leq{M}}\Re (\bar{C}_m)}{(\max\limits_{1\leq{m}\leq{M}}\vert\bar{C}_m\vert)^2}-\frac{\sqrt{2M_{max}}}{\pi\,d^{\frac{s}{t}}}\right)^{-2}\sum_{m=1}^{M}\left|Y_m\right|^2.
 \end{split}
\end{equation}

Thus, if $\frac{\min\limits_{1\leq{m}\leq{M}}\Re (\bar{C}_m)}{(\max\limits_{1\leq{m}\leq{M}}\vert\bar{C}_m\vert)^2}>\frac{\sqrt{2M_{max}}}{\pi\,d^{\frac{s}{t}}}$, then the matrix $\mathbf{B}$ 
in algebraic system \eqref{compacfrm1} is invertible.

\item
If $\Re (\lambda_{m, 0}) > 0$ for each $m$. In this case $\Re\bar{C}_m < 0$ and so to prove the invertibility of $\mathbf{B}$ in algebraic system \eqref{compacfrm1}, consider 
$(-\mathbf{C}^{-1}-\mathbf{B}_{n})X=Y$ in place of 
\eqref{systemsolve1-small-ac} and proceed in the way as it done for the case $\Re \lambda_m \geq 0$ for each $m$. Then we can get the following estimate% and then by following the way as it was done in \cite[Lemma 2.22]{DPC-SM:MMS2013}, we can prove that
 
\begin{equation}\label{fnlinvert-small-ac-2-f1}
\begin{split}
 \sum_{m=1}^{M}|X_m|^{2}
\leq 4\left(\frac{\min\limits_{1\leq{m}\leq{M}}\Re (-\bar{C}_m)}{(\max\limits_{1\leq{m}\leq{M}}\vert\bar{C}_m\vert)^2}-\frac{\sqrt{2M_{max}}}{\pi\,d^{\frac{s}{t}}}\right)^{-2}\sum_{m=1}^{M}\left|Y_m\right|^2.
 \end{split}
\end{equation}
and the invertibility of the matrix $\mathbf{B}$, under the assumption that $\frac{\min\limits_{1\leq{m}\leq{M}}\Re (-\bar{C}_m)}{(\max\limits_{1\leq{m}\leq{M}}\vert\bar{C}_m\vert)^2}>\frac{\sqrt{2M_{max}}}{\pi\,d^{\frac{s}{t}}}$.

\item We are left with the case where $\Re (\lambda_{m,0}) > 0$ for few $m$'s. In this case we multiply every line of the system (\ref{systemsolve1-small-ac}) corresponding to $\Re (\lambda_{m, 0}) > 0$ by $-1$. Hence, the system (\ref{systemsolve1-small-ac}) becomes
\begin{eqnarray}\label{systemsolve1-small-ac-1}(\tilde{\mathbf{C}}^{-1}+\tilde{\mathbf{B}}_{n})X&=&\tilde{Y}.\end{eqnarray}
where now every component of the diagonal matrix $\tilde{\mathbf{C}}^{-1}$ is positive and $\Vert \tilde{\mathbf{B}}_{n}\Vert_2 =\Vert \mathbf{B}_{n}\Vert_2 $. Then we are in the case (1).
\end{enumerate}
\end{proofe}

\begin{remark}\label{Remark-on-sign-lambda}
Assume that $\Re(\lambda_{m, 0}) <0$.  Following the computations in \cite[Lemma 2.22]{DPC-SM:MMS2013} and assuming that
$\frac{5\pi}{3}\frac{\min\limits_{1\leq{m}\leq{M}}\Re\bar{C}_m}{(\max\limits_{1\leq{m}\leq{M}}\vert\bar{C}_m\vert)^2}>\frac{\gamma}{d}$ and 
$\gamma:=\min\limits_{j\neq\,m,1\leq\,j,m\leq\,M}\cos(\kappa|z_m-z_j|) \geq 0$, we can prove that
%  \begin{equation}\label{fnlinvert-small-ac-2-f}
% \begin{split}
 $$\left(\frac{\min\limits_{1\leq{m}\leq{M}}\Re\bar{C}_m}{(\max\limits_{1\leq{m}\leq{M}}\vert\bar{C}_m\vert)^2}-\frac{3\gamma}{5\pi\,d}\right)\sum_{m=1}^{M}|X_m|^{2}
\leq2\left(\sum_{m=1}^{M}|X_m|^2\right)^{1/2}\left(\sum_{m=1}^{M}|Y_m|^2\right)^{1/2},
$$
% \end{split}
% \end{equation}
which yields
% \begin{equation}\label{fnlinvert-small-ac-2-f1}
% \begin{split}
$$
\sum_{m=1}^{M}|X_m|^{2}
\leq4 \left(\frac{\min\limits_{1\leq{m}\leq{M}}\Re\bar{C}_m}{(\max\limits_{1\leq{m}\leq{M}}\vert\bar{C}_m\vert)^2}-\frac{3\gamma}{5\pi\,d}\right)^{-2}\sum_{m=1}^{M}|Y_m|^2.
% \end{split}
$$
and thus the invertibility of the matrix $\mathbf{B}$ 
in the algebraic system \eqref{compacfrm1}. Observe that the condition 
$\frac{\min\limits_{1\leq{m}\leq{M}}\Re (\bar{C}_m)}{(\max\limits_{1\leq{m}\leq{M}}\vert\bar{C}_m\vert)^2}>\frac{\sqrt{2M_{max}}}{\pi\,d^{\frac{s}{t}}}$ 
is satisfied if $t<2-\beta$. 

We can see that if the number of the scatterers is $O(a^{-s})$ such that $s\leq 2-\beta$, then we do not need upper bound on $t$,
which allow as to have very close scatterers. However, if we want to have larger number of scatterers, i.e. $s$ limited only by the bound $s\leq 3$, then we need a condition on $t$, i.e.
$t\leq 2-\beta$, which makes more restrictions on the minimum distance between the scatterers.
\end{remark}

\bibliographystyle{abbrv}
% 
% \bibliography{Acoustic-Impedance-challa-references-LMF}

\end{document}